# Numerical Discretization Methods for Seismic Response Analysis of SDOF Systems: A Unified Perspective

Farid Ghahari[1,2*]

[1] *California Geological Survey, Sacramento, CA 95814, USA*

[2] *The B. John Garrick Institute for the Risk Sciences, University of California, Los Angeles, CA 90095, USA*

## ABSTRACT

This paper reviews the most commonly used numerical methods for solving the differential equation governing the dynamic response of linear elastic Single-Degree-of-Freedom (SDOF) systems. For more than 80 years since its introduction, the response spectrum—defined as the spectrum of the maximum response of a series of SDOF systems with different natural periods and/or damping ratios—has remained the cornerstone of every seismic design code. The second-order differential equation that governs the dynamic response of a linear elastic SDOF system must be solved numerically to generate such response spectra. Although only one or two well-accepted time-discretization methods have been predominantly used by the earthquake engineering community over the past decades, these methods are directly or indirectly related to a broader family of methods for solving Linear Time-Invariant (LTI) systems, which have been extensively applied in other branches of engineering, particularly electrical engineering. It has recently come to my attention that a portion of our community, particularly students, may not be fully familiar with these methods. In this paper, I review these methods and describe their mathematical background, with a focus on the relative displacement of the SDOF system under ground acceleration—an essential quantity for various types of response spectra. I also briefly review some of the numerical methods traditionally used within our community, highlighting their similarities and differences. I evaluate the accuracy of all numerical methods introduced in this paper through several examples with available analytical solutions. This study focuses on time-domain solutions that can be employed for real- or near-real-time response prediction, which is particularly important for applications such as earthquake early warning and post-earthquake assessment. The paper is written to enable readers to implement these methods with minimal effort; however, MATLAB codes for all methods discussed are also provided.

## INTRODUCTION

The Response Spectrum (RS), first introduced by Housner in 1941 (Housner, 1941), is a cornerstone of earthquake and structural engineering. It represents the maximum responses—such as relative displacement, absolute acceleration, or pseudo-acceleration—of a series of Single-Degree-of-Freedom (SDOF) systems with varying natural periods or frequencies and different damping levels. This tool quantifies the impact of earthquake ground motions on structural systems and has been integral to seismic design and assessment codes for decades (e.g., Chopra, 2007a). The RS is derived through time-history analysis of linear elastic SDOF systems, making it an essential method for evaluating structural performance under seismic loading.

**Figure 1** shows a schematic representation of a SDOF system with mass $m$, stiffness $k$, and damping $c$. The system is subjected to a horizontal ground acceleration $\ddot{u}_g(t)$, and its total displacement is denoted by $u^t(t)$. As shown in the figure, three forces act on the mass: the inertial force $m\ddot{u}^t(t)$, the damping force $c\dot{u}(t)$, and the elastic restoring strain force $ku(t)$. These forces must be in equilibrium at every instant in time. Therefore, the differential equation governing the behavior of the system can be expressed as

$$m\ddot{u}^t(t) + c\dot{u}(t) + ku(t) = 0. \tag{1}$$

---

[*] Research Scientist/Civil Engineer. Corresponding author: Farid Ghahari, E-mail: farid.ghahari@conservation.ca.gov

Knowing that $\ddot{u}^t(t) = \ddot{u}(t) + \ddot{u}_g(t)$ and dividing Eq. (1) by $m$, we obtain

$$\ddot{u}(t) + 2\xi\omega_n\dot{u}(t) + \omega_n^2 u(t) = -\ddot{u}_g(t), \tag{2}$$

where $\omega_n = \sqrt{k/m}$ is the natural frequency of the system in rad/s, and $\xi = c/2m\omega_n$ is the damping ratio. Eq. (2) is a fundamental equation in structural dynamics and earthquake engineering, appearing at the very beginning of most relevant textbooks. By solving this equation for a range of SDOFs with different $\omega_n$ but the same $\xi$ under a given ground-motion acceleration excitation, the response spectrum corresponding to that damping level can be constructed by plotting the maximum response of interest (e.g., relative displacement) versus the natural frequency, or more commonly, the natural period $T_n = 2\pi/\omega_n$ (see, e.g., Clough and Penzien, 1975).

The solution to the second-order differential equation above cannot be obtained in closed form unless the ground acceleration has an analytical representation, which is rarely the case. Therefore, Eq. (2) must be discretized in time and solved numerically, which inevitably introduces errors. Since there is no unique approach to time discretization, various methods have been developed, each with its advantages and limitations (see, e.g., Meirovitch, 1980). In this paper, I review some of the most commonly used methods within the earthquake and structural engineering community, while also drawing on the extensive body of knowledge from control theory and signal processing. Specifically, I show standard methods for discretizing continuous-time systems and derive the corresponding equations required to solve Eq. (2). From this general perspective, many of the methods currently used in earthquake engineering practice appear as special cases, which will be further discussed later. I also provide simple MATLAB (2024) codes that rely only on basic commands, without using built-in functions from specialized toolboxes.

It should be emphasized that this is not a comprehensive review of all numerical methods applicable to Eq. (2). For instance, frequency-domain methods are excluded because they cannot be applied to real-time or near-real-time analysis, as the complete ground-motion time history must be available. This limitation is particularly relevant given the growing interest in recent years in applications such as structural health monitoring (Ghahari et al., 2024a) and earthquake early warning (Ghahari et al., 2024b).

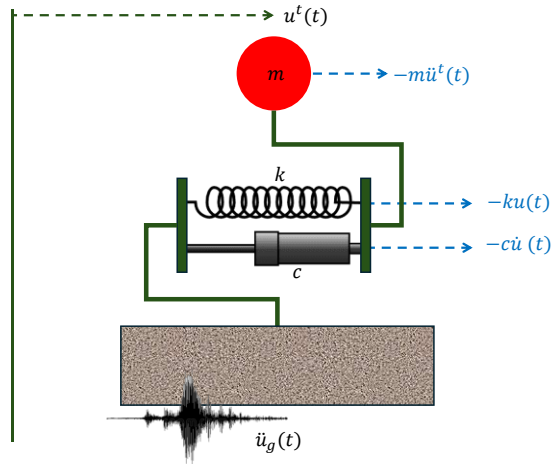

**Figure 1.** A schematic representation of an SDOF system

# TRANSFER FUNCTION APPROACH

A Transfer Function (TF) is a mathematical model that describes the relationship between the input and output of a Linear Time-Invariant (LTI) system. While the TF model has been employed in various fields, it is widely used in signal processing, communication theory, and control theory (see, e.g., Oppenheim et al., 1997). For the SDOF system introduced in Eq. (2), we can calculate the transfer function between the relative displacement and ground acceleration by applying the Laplace Transform on both sides of the equation. That is,

$$\ddot{U}(s) + 2\xi\omega_n \dot{U}(s) + \omega_n^2 U(s) = -\ddot{U}_g(s), \tag{3}$$

where capital letters represent the Laplace Transform, defined as:

$$F(s) = \int_0^\infty f(t)e^{-st}dt. \tag{4}$$

Assuming zero initial conditions, we have

$$\ddot{U}(s) = s^2 U(s), \dot{U}(s) = isU(s). \tag{5}$$

So, the transfer function can be calculated as

$$H(s) = \frac{U(s)}{\ddot{U}_g(s)} = \frac{-1}{s^2 + 2\xi\omega_n s + \omega_n^2}. \tag{6}$$

Having the Laplace Transform of the ground acceleration, this TF can be used to obtain the relative displacement time-history response of the system through inverse Laplace Transform. However, we only have access to a discrete-time version of the ground acceleration measured by digital accelerometers. In the discrete-time domain, there is a corresponding transfer function, which is represented through the z-Transform as:

$$H(z) = \frac{U(z)}{\ddot{U}_g(z)} = \frac{b_0 + b_1 z^{-1} + \cdots + b_{n_b} z^{-n_b}}{1 + a_1 z^{-1} + \cdots + a_{n_a} z^{-n_a}}, \tag{7}$$

where the z-Transform of a single-sided discrete-time signal $f[k]$ is defined as (see, e.g., Palani, 2021)

$$F(z) = \sum_{n=0}^\infty f[k]z^{-k}, \tag{8}$$

It is important to recall that any delay in the time domain appears as a power in the z-domain, i.e., $Z\{f[k - k_0]\} = z^{-k_0}F(z)$. So, provided that we have such a discrete-time transfer function, we can rewrite Eq. (7) as

$$(1 + a_1 z^{-1} + \cdots + a_{n_a} z^{-n_a})U(z) = (b_0 + b_1 z^{-1} + \cdots + b_{n_b} z^{-n_b})\ddot{U}_g(z), \tag{9}$$

And use the definition of the z-Transform and the mentioned property to write the following equation:

$$u[k] = -a_1 u[k-1] - \cdots - a_{n_a} u[k - n_a] + b_0 \ddot{u}_g[k] + b_1 \ddot{u}_g[k-1] + \cdots + b_{n_b} \ddot{u}_g[k - n_b], \tag{10}$$

which is a recursive formulation to calculate the relative displacement of the SDOF in real time. However, contrary to the continuous-time transfer function (Eq. (6)), the discrete-time transfer function (Eq. (7)) cannot be directly deduced from the governing equation. To obtain a discrete-time transfer function, either the equation of motion or the continuous-time transfer function has to be discretized in time. In what follows, several of such discretization methods will be presented.

*Zero-Order Hold Method*

As mentioned earlier, having the Laplace Transform of the input excitation, the response can be easily computed as $L^{-1}\{H(s)\ddot{U}_g(s)\}$, where $L^{-1}\{.\}$ stands for the inverse Laplace Transform. The response would be identical to the transfer function if the $\ddot{U}_g(s) = 1$, which is the case for $\ddot{u}_g(t) = \delta(t)$ where

$$\delta(t) = \begin{cases} 0, & t \neq 0 \\ \infty, & t = 0 \end{cases} \quad (11)$$

with $\int_{-\infty}^{\infty} \delta(t)dt = 1$. So, if we apply a discrete-time impulse whose reconstructed continuous-time is $\delta(t)$, the sampled version of the output signal would be the discrete-time transfer function. One method to convert a discrete-time signal to a continuous-time version is by holding each sample value for one sample interval, which is called Zero-Order-Hold (ZOH). Through this method, a discrete-time signal $f[k]$ sampled at time instants $k\Delta t$ is converted to a continuous-time signal $f(t)$ as

$$f(t) = \sum_{n=0}^{\infty} f[k]\, rect\left(\frac{t - \frac{\Delta t}{2} - k\Delta t}{\Delta t}\right), \quad (12)$$

where $rect\left(\frac{t - \frac{\Delta t}{2}}{\Delta t}\right)$ is a function with amplitude 1 from 0 to $\Delta t$ and zero elsewhere. Assuming $f[k] = \delta[k]$ (discrete-time impulse) with

$$\delta[k] = \begin{cases} 0, & k \neq 0 \\ 1, & k = 0 \end{cases} \quad (13)$$

we have

$$L\left\{\sum_{n=0}^{\infty} \delta[k]\, rect\left(\frac{t - \frac{\Delta t}{2} - k\Delta t}{\Delta t}\right)\right\} = \frac{1 - e^{-s\Delta t}}{s}. \quad (14)$$

Therefore, the discrete-time TF can be calculated as

$$H(z) = (1 - z^{-1})\, Z\left\{L^{-1}\left\{\frac{H(s)}{s}\right\}\bigg|_{t=k\Delta t}\right\}, \quad (15)$$

which can be represented as a rational function like Eq. (7) with the coefficients reported in Table 1, where $\omega_d = \omega_n\sqrt{1 - \xi^2}$.

Table 1. Coefficients of the discrete-time transfer function discretized using the ZOH method.

| Coefficient | Value |
|---|---|
| $a_1$ | $1 - 2e^{-\xi\omega_n\Delta t}\cos\omega_d\Delta t$ |
| $a_2$ | $e^{-2\xi\omega_n\Delta t}$ |
| $b_0$ | 0 |
| $b_1$ | $-\frac{1}{\omega_n^2}(1 - e^{-\xi\omega_n\Delta t}\cos\omega_d\Delta t - \frac{\xi}{\sqrt{1-\xi^2}}e^{-\xi\omega_n\Delta t}\sin\omega_d\Delta t)$ |
| $b_2$ | $-\frac{1}{\omega_n^2}(e^{-2\xi\omega_n\Delta t} - e^{-\xi\omega_n\Delta t}\cos\omega_d\Delta t + \frac{\xi}{\sqrt{1-\xi^2}}e^{-\xi\omega_n\Delta t}\sin\omega_d\Delta t)$ |

*First-Order Hold Method*

We can follow a similar approach to convert a discrete-time impulse to a continuous impulse, but assuming the signal varies linearly between every two consecutive samples. This type of signal reconstruction is called the First-Order-Hold (FOH) reconstruction. In such a case, Eq. (15) will change to

$$H(z) = \frac{(1-z^{-1})^2}{\Delta t} Z\left\{L^{-1}\left\{\frac{H(s)}{s^2}\right\}\bigg|_{t=k\Delta t}\right\}, \quad (16)$$

which can be represented as a rational function like Eq. (7) with the coefficients reported in Table 2.

Table 2. Coefficients of the discrete-time transfer function discretized using the FOH method.

| Coefficient | Value |
|---|---|
| $a_1$ | $1 - 2e^{-\xi\omega_n\Delta t}\cos\omega_d\Delta t$ |
| $a_2$ | $e^{-2\xi\omega_n\Delta t}$ |
| $b_0$ | $\frac{1}{\Delta t}\left[\frac{2\xi}{\omega_n^3}\left(1 - e^{-\xi\omega_n\Delta t}\cos\omega_d\Delta t\right) - \frac{\Delta t}{\omega_n^2} + \frac{1-2\xi^2}{\omega_n^2\omega_d}e^{-\xi\omega_n\Delta t}\sin\omega_d\Delta t\right]$ |
| $b_1$ | $\frac{2}{\Delta t}\left[\frac{\xi}{\omega_n^3}\left(e^{-2\xi\omega_n\Delta t} - 1\right) + \frac{\Delta t}{\omega_n^2}e^{-\xi\omega_n\Delta t}\cos\omega_d\Delta t - \frac{1-2\xi^2}{\omega_n^2\omega_d}e^{-\xi\omega_n\Delta t}\sin\omega_d\Delta t\right]$ |
| $b_2$ | $\frac{1}{\Delta t}\left[-\left(\frac{2\xi}{\omega_n^3} + \frac{\Delta t}{\omega_n^2}\right)e^{-2\xi\omega_n\Delta t} + \frac{2\xi}{\omega_n^3}e^{-\xi\omega_n\Delta t}\cos\omega_d\Delta t + \frac{1-2\xi^2}{\omega_n^2\omega_d}e^{-\xi\omega_n\Delta t}\sin\omega_d\Delta t\right]$ |

*Impulse-Invariant Method*

It is a well-proven fact that multiplication in the Laplace domain is equivalent to the convolution in the time domain. So, we can write

$$L^{-1}\{H(s)\ddot{U}_g(s)\} = L^{-1}\{H(s)\}(t) * \ddot{u}_g(t). \quad (17)$$

The linear convolution on the right-hand side of Eq. (17) is called the Duhamel integral (Chopra, 2007b), and the inverse Laplace Transform of the transfer function of the SDOF system is called the Impulse Response Function (IRF), which can be analytically calculated as

$$L^{-1}\{H(s)\}(t) = \frac{-1}{\omega_d}e^{-\xi\omega_n t}\sin\omega_d t. \quad (18)$$

One way to discretize the continuous-time transfer function is to make sure the discrete-time IRF exactly matches the continuous-time IRF at the samples. The discrete-time IRF can be written as

$$h[k] = \frac{-\Delta t}{\omega_d}e^{-\xi\omega_n k\Delta t}\sin\omega_d k\Delta t. \quad (19)$$

By using the z-Transform, it is easy to show that the discrete-time transfer function would be a second-order rational function as Eq. (7) with the coefficients reported in Table 3.

Table 3. Coefficients of the discrete-time transfer function discretized using the Impulse-Invariant method.

| Coefficient | Value |
|---|---|
| $a_1$ | $-2e^{-\xi\omega_n\Delta t}\cos\omega_d\Delta t$ |
| $a_2$ | $e^{-2\xi\omega_n\Delta t}$ |
| $b_0$ | 0 |
| $b_1$ | $-\frac{\Delta t}{\omega_d}e^{-\xi\omega_n\Delta t}\sin\omega_d\Delta t$ |
| $b_2$ | 0 |

*Euler Method*

Having a continuous-time transfer function, one idea would be to directly use $z = e^{s\Delta t}$ to discretize it. However, this discretization would result in a nonlinear numerator and denominator. To resolve this issue, we can use the Taylor expansion:

$$e^{s\Delta t} = 1 + s\Delta t + \frac{(s\Delta t)^2}{2!} + \cdots \tag{20}$$

and approximate $e^{s\Delta t}$ with the first two terms to obtain $s \approx \frac{z-1}{\Delta t}$. Replacing this relationship in Eq. (6), we find a discrete-time transfer function as represented in Eq. (7) with the coefficients reported in Table 4. This method is called the Forward Euler method.

Table 4. Coefficients of the discrete-time transfer function discretized using the Forward Euler method.

| Coefficient | Value |
| --- | --- |
| $a_1$ | $2\xi\omega_n\Delta t - 2$ |
| $a_2$ | $1 - 2\xi\omega_n\Delta t + \omega_n^2\Delta t^2$ |
| $b_0$ | 0 |
| $b_1$ | 0 |
| $b_2$ | $-\Delta t^2$ |

It is also easy to show that the $\frac{1}{1-s\Delta t}$ has a Taylor expansion with exactly the same first two terms as $e^{s\Delta t}$. So, we can also use $z \approx \frac{1}{1-s\Delta t}$ to discretize the transfer function, which is called the Backward Euler method. The coefficients of this transfer function are reported in Table 5.

Table 5. Coefficients of the discrete-time transfer function discretized using the Backward Euler method. In this table, $\rho = 1 + 2\xi\omega_n\Delta t + \omega_n^2\Delta t^2$.

| Coefficient | Value |
| --- | --- |
| $a_1$ | $-2(1 + \xi\omega_n\Delta t)/\rho$ |
| $a_2$ | $1/\rho$ |
| $b_0$ | $-\Delta t^2/\rho$ |
| $b_1$ | 0 |
| $b_2$ | 0 |

*Tustin Method*

To increase the accuracy of the approximation in the Euler method, we can write

$$z = e^{s\Delta t} = \frac{e^{s\frac{\Delta t}{2}}}{e^{-s\frac{\Delta t}{2}}}. \tag{21}$$

Then, we can use the first two terms of the Taylor expansion of the numerator and denominator, that is

$$z \approx \frac{1+\frac{s\Delta t}{2}}{1-\frac{s\Delta t}{2}}, \tag{22}$$

from which we have $s \approx \frac{2}{\Delta t}\frac{z-1}{z+1}$. By using this relationship, the coefficients of the transfer function will be those reported in Table 6. This method is called the Tustin (Tustin, 1947), bilinear, or trapezoidal method.

Table 6. Coefficients of the discrete-time transfer function discretized using the Tustin method. In this table, $\rho = 4 + 4\xi\omega_n\Delta t + \omega_n^2\Delta t^2$.

| Coefficient | Value |
|---|---|
| $a_1$ | $(2\omega_n^2\Delta t^2 - 8)/\rho$ |
| $a_2$ | $(4 - 4\xi\omega_n\Delta t + \omega_n^2\Delta t^2)/\rho$ |
| $b_0$ | $-\Delta t^2/\rho$ |
| $b_1$ | $-2\Delta t^2/\rho$ |
| $b_2$ | $-\Delta t^2/\rho$ |

When using the Tustin transformation, the entire imaginary axes (frequency axes) of the Laplace domain are compressed into the $2\pi$-length of the unit circle in the z-domain, causing a frequency distortion. Therefore, if there is an important frequency (in our case, the natural frequency of the SDOF system at which the continuous-time transfer function shows a peak), the location of that frequency might move in the discrete-time transfer function, and such a distortion would be larger if that frequency is closer to the Nyquist frequency. To mathematically show this phenomenon, let's assume that the frequency of interest is $\omega_n$. Replacing $s = i\omega_n$ into Eq. (22), the corresponding frequency in the discrete-time domain can be calculated as

$$e^{i\widetilde{\omega}_n\Delta t} \approx \frac{1+\frac{i\omega_n\Delta t}{2}}{1-\frac{i\omega_n\Delta t}{2}}, \tag{23}$$

which results in

$$\widetilde{\omega}_n \approx \frac{1}{i\Delta t}\ln\frac{1+\frac{i\omega_n\Delta t}{2}}{1-\frac{i\omega_n\Delta t}{2}}, \tag{24}$$

which is different from $\omega_n$ unless $\Delta t$ is very small. To prevent such distortion at this frequency of interest, we can pre-warp the transformation as $s \approx \frac{\omega_n}{\tan\left(\frac{\omega_n\Delta t}{2}\right)}\frac{z-1}{z+1}$. The corresponding coefficients of the discrete-time transfer function are reported in Table 7.

Table 7. Coefficients of the discrete-time transfer function discretized using the Tustin method with pre-warping. In this table, $\eta = \frac{\omega_n}{\tan\left(\frac{\omega_n \Delta t}{2}\right)}$, and $\rho = \eta^2 + 2\eta\xi\omega_n + \omega_n^2$.

| Coefficient | Value |
|---|---|
| $a_1$ | $(2\omega_n^2 - 2\eta^2)/\rho$ |
| $a_2$ | $(\eta^2 - 2\eta\xi\omega_n + \omega_n^2)/\rho$ |
| $b_0$ | $-1/\rho$ |
| $b_1$ | $-2/\rho$ |
| $b_2$ | $-1/\rho$ |

*Matched Method*

The dynamic of LTI systems is defined by the locations of their poles (roots of denominator), zeros (roots of numerator), and DC gain (transfer function value at $s = 0$). For our continuous-time SDOF transfer function, the system does not have any finite zeros, has two poles at $-\xi\omega_n \mp i\omega_d$, and a DC gain of $-1/\omega_n^2$. We can determine the exact locations of these zeros and poles in the discrete-time domain using the nonlinear mapping of $z = e^{s\Delta t}$. Therefore, the transformed poles will be $p_{1,2} = e^{(-\xi\omega_n \mp i\omega_d)\Delta t}$. While the continuous-time system does not have finite zeros, since the difference between the orders of the numerator and denominator is 2, it indeed has two zeros at infinity. In the z-domain, these zeros are mapped to -1 (replace $s = \pm\infty$ into Eq. (22)), that is, $z_{1,2} = -1$. So, the initial representation of the discrete-time transfer function is

$$H(z) = \frac{(1-z_1 z^{-1})(1-z_2 z^{-1})}{(1-p_1 z^{-1})(1-p_2 z^{-1})}. \tag{25}$$

To make sure the system is causal and strictly proper (Oppenheim, 1999), the number of zeros must be less than the number of poles. Therefore, we remove one of the infinite zeros, resulting in

$$H(z) = \frac{(1-z_1 z^{-1})}{(1-p_1 z^{-1})(1-p_2 z^{-1})}. \tag{26}$$

Finally, the gain must match that of the continuous-time transfer function. The DC gain of the discrete-time transfer function in Eq. (26) is $H(z=1) = 2/(1-p_1 z^{-1})(1-p_2 z^{-1})$, so we multiply the transfer function by $K = -(1-p_1 z^{-1})(1-p_2 z^{-1})/2\omega_n^2$. The coefficients of the final discrete-time transfer function are reported in Table 8.

Table 8. Coefficients of the discrete-time transfer function discretized using the Matched method.

| Coefficient | Value |
|---|---|
| $a_1$ | $-2e^{-\xi\omega_n\Delta t} \cos \omega_d \Delta t$ |
| $a_2$ | $e^{-2\xi\omega_n\Delta t}$ |
| $b_0$ | 0 |
| $b_1$ | $\left(2e^{-\xi\omega_n\Delta t} \cos \omega_d \Delta t - e^{-2\xi\omega_n\Delta t} - 1\right)/(2\omega_n^2)$ |
| $b_2$ | $\left(2e^{-\xi\omega_n\Delta t} \cos \omega_d \Delta t - e^{-2\xi\omega_n\Delta t} - 1\right)/(2\omega_n^2)$ |

*Least-Squares Method*

The least squares method minimizes the error between the frequency responses of the continuous-time and discrete-time systems up to the Nyquist frequency using a vector-fitting optimization approach. This method is useful when you want to capture fast system dynamics but must use a larger sample time. When using the same sample time as the Tustin approximation or Matched method, you get a smaller difference between the continuous-time and discrete-time frequency responses through this approach.

To use this method, the Frequency Response Function (FRF) of the continuous-time transfer function is calculated over a range of frequencies of interest by replacing $s = i\omega$ into Eq. (6):

$$F(\omega) = \frac{-1}{\omega_n^2 - \omega^2 + 2i\xi\omega_n}. \tag{27}$$

Then, the corresponding discrete-time FRF values are calculated by replacing $z = e^{i\omega\Delta t}$ into Eq. (7):

$$F_D(\omega) = \frac{B(\omega)}{A(\omega)}, \tag{28}$$

where

$$A(\omega) = 1 + a_1 e^{-i\omega\Delta t} + a_2 e^{-2i\omega\Delta t}, \tag{29}$$

$$B(\omega) = b_0 + b_1 e^{-i\omega\Delta t} + b_2 e^{-2i\omega\Delta t}. \tag{30}$$

Finally, we solve the following minimization to obtain unknown coefficients:

$$\min_{\substack{a_1, a_2, \\ b_0, b_1, b_2}} \sum_q W(\omega_q) \left| \frac{B(\omega_q)}{A(\omega_q)} - F(\omega_q) \right|^2, \tag{31}$$

where $W(\omega_q)$ is the weight of the $qth$ frequency. This nonlinear minimization can be converted to a linear minimization by internally changing the weight of each term and defining an equivalent minimization as below

$$\min_{\substack{a_1, a_2, \\ b_0, b_1, b_2}} \sum_q W(\omega_q) |B(\omega_q) - A(\omega_q) F(\omega_q)|^2, \tag{32}$$

which can be represented in a vector form as

$$\min_{\theta} \|\mathbf{D}\boldsymbol{\theta} - \mathbf{h}\|^2, \tag{33}$$

where $\boldsymbol{\theta} = [a_1 \ a_2 \ b_0 \ b_1 \ b_2]^T$ and matrix $\mathbf{D}$ and vector $\mathbf{h}$ are constructed using $z_q = e^{i\omega_q \Delta t}$, $F_q = F(\omega_q)$, and $W_q = \sqrt{W(\omega_q)}$ terms as below

$$\mathbf{D} = \begin{bmatrix} W_1 F_1 z_1^{-1} & W_1 F_1 z_1^{-2} & W_1 & W_1 z_1^{-1} & W_1 z_1^{-2} \\ \vdots & \vdots & \vdots & \vdots & \vdots \\ W_N F_N z_N^{-1} & W_N F_N z_N^{-2} & W_N & W_N z_N^{-1} & W_N z_N^{-2} \end{bmatrix}, \tag{34}$$

$$\mathbf{h} = \begin{bmatrix} W_1 F_1 \\ \vdots \\ W_N F_N \end{bmatrix}. \tag{35}$$

For real-valued $\mathbf{D}$ and $\mathbf{h}$, the solution to Eq. (33) can be calculated as $\boldsymbol{\theta} = (\mathbf{D}^T \mathbf{D})^{-1} \mathbf{D}^T \mathbf{h}$. In our case, $\mathbf{D}$ and $\mathbf{h}$ are complex-valued matrices and vectors, respectively, but can be expanded into real and imaginary parts and solved as described in (Levi, 1959). However, because the frequency response function is complex conjugate symmetric, we can use only positive frequencies and solve it as

$$\boldsymbol{\theta} = (Real\{\mathbf{D}^T \mathbf{D}\})^{-1} Real\{\mathbf{D}^T \mathbf{h}\}. \tag{36}$$

The stability of all the discretization methods will be discussed later, and this method is among the ones that does not inherently preserve the stability of the continuous-time system. In other words, the stability condition—that the poles of a discrete-time system must be inside the unit circle—could be violated. To resolve this issue, the minimization in Eq. (31) is solved through a nonlinear optimization using the damped Gauss-Newton method (Dennis and Schnabel, 1983) during which any poles $p_i$ outside the unit circle are reflected back into the unit circle as $p_i \approx p_i/|p_i|^2$. To start the minimization process, solutions obtained from Eq. (36) can be employed.

While the MATLAB code is provided in this paper to construct the discrete-time transfer function for any system, Table 9 presents the transfer function coefficients for a range of common natural periods with $\xi = 5\%$, and $\Delta t = 0.01\ sec$.

Table 9. Coefficients of the discrete-time transfer function discretized using the Least-Squares method.

| $T_n(sec.)$ | $a_1$ | $a_2$ | $b_0\ (\times 10^4)$ | $b_1\ (\times 10^4)$ | $b_2\ (\times 10^4)$ |
|---|---|---|---|---|---|
| 0.05 | -0.5830 | 0.8825 | -0.1072 | -0.6442 | -0.0976 |
| 0.075 | -1.2848 | 0.9197 | -0.1038 | -0.7147 | -0.0980 |
| 0.1 | -1.5689 | 0.9391 | -0.1027 | -0.7435 | -0.0985 |
| 0.15 | -1.7896 | 0.9590 | -0.1019 | -0.7671 | -0.0992 |
| 0.2 | -1.8727 | 0.9691 | -0.1017 | -0.7769 | -0.0996 |
| 0.25 | -1.9131 | 0.9752 | -0.1015 | -0.7820 | -0.0999 |
| 0.3 | -1.9360 | 0.9793 | -0.1015 | -0.7852 | -0.1001 |
| 0.4 | -1.9600 | 0.9844 | -0.1014 | -0.7888 | -0.1004 |
| 0.5 | -1.9718 | 0.9875 | -0.1014 | -0.7908 | -0.1006 |
| 0.75 | -1.9847 | 0.9917 | -0.1013 | -0.7932 | -0.1008 |
| 1 | -1.9898 | 0.9937 | -0.1013 | -0.7943 | -0.1009 |
| 1.5 | -1.9941 | 0.9958 | -0.1013 | -0.7954 | -0.1011 |
| 2 | -1.9959 | 0.9969 | -0.1013 | -0.7959 | -0.1011 |
| 3 | -1.9975 | 0.9979 | -0.1013 | -0.7964 | -0.1012 |
| 4 | -1.9982 | 0.9984 | -0.1013 | -0.7967 | -0.1012 |
| 5 | -1.9986 | 0.9987 | -0.1013 | -0.7968 | -0.1012 |
| 7.5 | -1.9991 | 0.9992 | -0.1013 | -0.7970 | -0.1013 |
| 10 | -1.9993 | 0.9994 | -0.1013 | -0.7971 | -0.1013 |

*Kanamori et al. (1999) Method*

The last discretization method presented under the transfer function family is a special case of the previous least-squares method. Kanamori et al. (1999) developed a recursive filter similar to the one presented in (10) to calculate the displacement response of an SDOF system with natural frequency $\omega_n$ and damping ratio $\xi$. This method which has been widely adopted in strong motion community (see, e.g., Kohler et al., 2020), is based on a discrete-time transfer function whose coefficients are calculated by minimizing the difference between the analytical continuous-time transfer function and a parametric discrete-time transfer function with two main differences: 1- only amplitude of the FRF is used, and 2- instead of the generic representation of the discrete-time filter shown in Eq. (7), Kanamori et al. (1999) derived the transfer function parametric representation by discretizing differential equation of motion (Eq. (2)) using the Backward Finite Difference (BFD) method. Eq. (2) discretized in time using BFD can be written as

$$\left(\frac{u[k]-2u[k-1]+u[k-2]}{\Delta t^2}\right) + 2\xi\omega_n\left(\frac{u[k]-u[k-1]}{\Delta t}\right) + \omega_n^2 u[k] = -\ddot{u}_g[k]. \tag{37}$$

Applying the z-Transform, we have

$$[U(z) - 2z^{-1}U(z) + z^{-2}U(z)] + 2\xi\omega_n\Delta t[U(z) - z^{-1}U(z)] + \omega_n^2\Delta t^2 U(z) = -\Delta t^2 \ddot{U}_g(z), \tag{38}$$

from which the discrete-time transfer function can be expressed as

$$H(z) = \frac{-\Delta t^2}{1 + 2\xi\omega_n\Delta t + \omega_n^2\Delta t^2 - 2(1+\xi\omega_n\Delta t)z^{-1} + z^{-2}}. \tag{39}$$

There are no unknown coefficients in Eq. (39), but as this is an approximation of the continuous-time transfer function, Kanamori et al. (1999) estimated equivalent natural frequency and damping ratio by solving the following minimization

$$\min_{\tilde{\xi},\tilde{\omega}_n} \left\| \left|\frac{-1}{-\omega^2 + 2i\xi\omega_n\omega + \omega_n^2}\right| - \left|\frac{-\Delta t^2}{1 + 2\tilde{\xi}\tilde{\omega}_n\Delta t + \tilde{\omega}_n^2\Delta t^2 - 2(1+\tilde{\xi}\tilde{\omega}_n\Delta t)e^{-i\omega\Delta t} + e^{-2i\omega\Delta t}}\right| \right\|, \tag{40}$$

for $0.01 \leq \frac{\omega}{2\pi} \leq 10$, and reported the results for three cases using a sampling rate of 100 Hz, as shown in Table 10. The reported error across the specified frequency range is less than 3%.

Table 10. Effective natural frequencies and damping ratios obtained using Kanamori et al. (1999)

| Actual SDOF | | Effective SDOF | |
|---|---|---|---|
| $\omega_n$ (rad/sec.) | $\xi$ (%) | $\tilde{\omega}_n$ (rad/sec.) | $\tilde{\xi}$ (%) |
| 20.94 | 5.00 | 21.01 | -5.41 |
| 6.28 | 5.00 | 6.28 | 1.80 |
| 2.09 | 5.00 | 2.10 | 4.40 |

A MATLAB function implementing all the discretization methods discussed in this section is presented in **Figure 2**. After discretizing the transfer function using one of these methods, the function computes the relative displacement of the structure. While additional pre- and post-processing steps will be discussed later in the "Implementation Remarks" section, two separate MATLAB functions are provided in **Figure 3** and **Figure 4**, corresponding to the least-squares and Kanamori et al. (1999) methods, respectively. The function in **Figure 3** implements the solution of Eq. (33) using a nonlinear optimization approach in which the stability of the transfer function is enforced, whereas the function in **Figure 4** evaluates the objective function defined in Eq. (40).

```matlab
function [disp,Time] = TransferFunction(wn,xi,ag,dt,dta,dtd,c2d_method,interp_flag)
% Inputs:
%           % wn : natural frequency (rad/s)
%           % xi : damping ratio
%           % ag : ground acceleration (Nx1)
%           % dt : ground motion time step (s)
%           % dta : analysis time step (s), used if interp_flag = 1 or 2
%           % dtd : display time step (s)
%           % c2d_method : method used for time discretization
%           % interp_flag : interpolation flag (0: no interpolation, 1: linear, 2: sinc)

% Outputs:
%           % disp : relative displacement response (Nx1)
%           % Time : time vector (Nx1)

if interp_flag==1
        US = dt/dta;
        N = length(ag);
        Time = [0:dt:(N-1)*dt]';
        ta = [0:dta:(US*N-1)*dta]';
        ag = interp1(Time,ag,ta,'linear', 'extrap');
        dt = dta;
elseif interp_flag==2
        US = dt/dta;
        ag = interp(ag,US);
        dt = dta;
end
N = length(ag);
Time = [0:dt:(N-1)*dt]';
alpha = exp(-xi*wn*dt);
theta = wn*sqrt(1-xi^2)*dt;
b = xi/sqrt(1-xi^2);
wd = wn*sqrt(1-xi^2);
if strcmp(c2d_method,'ZOH')
        Num = -1/(wn^2)*[0,1-alpha*cos(theta)-b*alpha*sin(theta),alpha^2-alpha*cos(theta)+b*alpha*sin(theta)];
        Den = [1,-2*alpha*cos(theta),alpha^2];
elseif strcmp(c2d_method,'FOH')
        b0 = 1/dt*(2*xi/wn^3*(1-alpha*cos(theta))-dt/wn^2+(1-2*xi^2)/(wn^2*wd)*alpha*sin(theta));
        b1 = 1/dt*(2*xi/wn^3*(alpha^2-1)+2*dt/wn^2*alpha*cos(theta)-2*(1-2*xi^2)/(wn^2*wd)*alpha*sin(theta));
        b2 = 1/dt*(-(2*xi/wn^3+dt/wn^2)*alpha^2+2*xi/wn^3*alpha*cos(theta)+(1-2*xi^2)/(wn^2*wd)*alpha*sin(theta));
        Num = [b0,b1,b2];
        Den = [1,-2*alpha*cos(theta),alpha^2];
elseif strcmp(c2d_method,'FE')
        Num = [0,0,-dt^2];
        Den = [1,2*xi*wn*dt-2,1-2*xi*wn*dt+wn^2*dt^2];
elseif strcmp(c2d_method,'BE')
        Num = 1/(1+2*xi*wn*dt+wn^2*dt^2)*[0,0,-dt^2];
        Den = [1,-(2+2*xi*wn*dt)/(1+2*xi*wn*dt+wn^2*dt^2),1/(1+2*xi*wn*dt+wn^2*dt^2)];
elseif strcmp(c2d_method,'Matched')
        Num = -1/(2*wn^2)*(1-2*exp(-xi*wn*dt)*cos(wd*dt)+exp(-2*xi*wn*dt))*[0,1,1];
        Den = [1,-2*exp(-xi*wn*dt)*cos(wd*dt),exp(-2*xi*wn*dt)];
elseif strcmp(c2d_method,'Tustin')
        Num = -dt^2*[1,2,1]/(4+4*dt*xi*wn+wn^2*dt^2);
        Den = [1,(2*wn^2*dt^2-8)/(4+4*dt*xi*wn+wn^2*dt^2),(4-4*dt*xi*wn+wn^2*dt^2)/(4+4*dt*xi*wn+wn^2*dt^2)];
elseif strcmp(c2d_method,'Tustin2')
        eta = wn/(tan(wn*dt/2));
        Ro = (eta^2+2*xi*wn*eta+wn^2);
        Num = -[1,2,1]/Ro;
        Den = [1,(2*wn^2-2*eta^2)/Ro,(wn^2+eta^2-2*xi*wn*eta)/Ro];
elseif strcmp(c2d_method,'Impulse')
        Num = -1/(wd)*[0,alpha*dt*sin(theta),0];
        Den = [1,-2*alpha*cos(theta),alpha^2];
elseif strcmp(c2d_method,'LSQ')
        w = 2*pi*[0:0.01:(1/dt)/2]';
        s = 1i*w;
        h = -1./(s.^2+2*xi*wn.*s+wn^2);
        W = ones(length(w),1);
        [Num,Den] = My_invfreqz(h,w/(1/dt),2,2,W,100,1e-5);
elseif strcmp(c2d_method,'Kanamori')
        w = 2*pi*[0:0.01:25]';
        s = 1i*w;
        hc = -1./(s.^2+2*xi*wn.*s+wn^2);
        options = optimset('TolFun',1e-10,'TolX',1e-6,'MaxFunEvals',1000,'MaxIter',1000);
        x = fmincon(@(x)Kanamori(x,w,hc,dt),[wn;xi],[-dt/2-1],0,[],[],[0.5*wn,-1000],[1.5*wn,0.99],[],options);
        w0 = x(1); h = x(2);
        Num = -[dt^2,0,0]/(1+2*h*w0*dt+w0^2*dt^2);
        Den = [1,-2*(1+h*w0*dt)/(1+2*h*w0*dt+w0^2*dt^2),1/(1+2*h*w0*dt+w0^2*dt^2)];
elseif strcmp(c2d_method,'CD')
        Num = [0,-dt^2/(1+xi*wn*dt),0];
        Den = [1,(wn^2*dt^2-2)/(1+xi*wn*dt),(1-xi*wn*dt)/(1+xi*wn*dt)];
end

disp = zeros(N,1);
for k=3:N
        disp(k,1) = -Den(2)*disp(k-1)-Den(3)*disp(k-2)+Num(1)*ag(k)+Num(2)*ag(k-1)+Num(3)*ag(k-2);
end

if dt>dtd
        US = dt/dtd;
        disp = interp(disp,US);
        Time = [0:dtd:(US*N-1)*dtd]';
end
end
```

**Figure 2.** Transfer function MATLAB function for calculating the time-history response of an SDOF system.

```matlab
function [b,a] = My_invfreqz(H,w,nb,na,W,maxiter,tol)
% This function calculates the real-valued coefficients of a discrete-time
% filter fitted to the frequency response function data points given in H.
% This code is adapted from the built-in "invfreqz" function in MATLAB.
% Inputs:
        % w: Normalized positive discrete frequencies (Lx1)
        % H: Impedance function at given discrete frequencies (Lx1)
        % nb: Order of the numerator
        % na: Order of the denominator
        % W: Weight function (Lx1)
        % maxiter: Maximum number of iteration for the stability enforcement
        % tol: Error tolearnce for the stability enforcement
% Outputs:
        % b: Numerator coefficients of the estimated discrete-time filter (nb+1x1)
        % a: Denominator coefficients of the estimated discrete-time filter (nax1)
if isempty(W)
                W = ones(length(w),1);
end
Wf = diag(sqrt(W));
nm = max(na,nb+1);
Z = exp(-1i*(0:nm)'.'*w');
Da = -(Z(2:na+1,:).').*(H*ones(1,na));
Db = (Z(1:nb+1,:).');
D = Wf*[Da Db];
R = real(D'*D);
Vd = real(D'*(Wf*H));
Th = R\Vd;
a = [1 th(1:na).'];
b = [th(na+1:na+nb+1).'];
if nargin>5 % iterative process
        indb = 1:length(b);
        indg = 1:length(a);
        % Stabilizing the denominator
        a = polystab(a);
        % The initial estimate:
        GC = ((b*Z(indb,:))./(a*Z(indg,:))).';
        e = Wf*(GC-H);
        MaxError = e'*e;
        th = [a(2:na+1) b(1:nb+1)].';
        % Newton-Gauss minimization
        dth = 2*tol+1;
        l = 0;
        st = 0;
        while all([norm(dth)>tol l<maxiter st~=1])
                DH_a =(Z(2:na+1,:).').*(-GC./((a*Z(1:na+1,:)).')*ones(1,na));
                DH_b =(Z(1:nb+1,:).')./((a*Z(1:na+1,:)).'*ones(1,nb+1));
                J = Wf*[DH_a,DH_b];
                e = Wf*(GC-H);
                R = real(J'*J);
                Vd = real(J'*e);
                dth = R\Vd;
                ll = 0;
                k = 1;
                V1 = MaxError+1;
                while all([V1 >= MaxError ll<20])
                        t1 = th-k*dth;
                        if ll==19
                                t1=th;
                        end
                        a = polystab([1 t1(1:na).']);
                        t1(1:na) = a(2:na+1).';
                        b = t1(na+1:na+nb+1).';
                        GC = ((b*Z(indb,:))./(a*Z(indg,:))).';
                        e = Wf*(GC-H);
                        V1 = e'*e;
                        t1 = [a(2:na+1) b(1:nb+1)].';
                        k = k/2;
                        ll=ll+1;
                        if ll==20
                                st=1;
                        end
                        if ll==10
                                dth = Vd/norm(R)*length(R);
                                k=1;
                        end
                end
                th = t1;
                MaxError = V1;
        end
end
```

**Figure 3.** Function used within least-squares method.

```
function residual = Kanamori(x,w,hc,dt)
        wn = x(1);
        xi = x(2);
        Num = -[dt^2,0,0]/(1+2*xi*wn*dt+wn^2*dt^2);
        Den = [1,-2*(1+xi*wn*dt)/(1+2*xi*wn*dt+wn^2*dt^2),1/(1+2*xi*wn*dt+wn^2*dt^2)];
        z = exp(1i.*w.*dt);
        hd = (Num(1)*z.^(0)+Num(2)*z.^(-1)+Num(3)*z.^(-2))./(Den(1)*z.^(0)+Den(2)*z.^(-1)+Den(3)*z.^(-2));
        residual = norm(abs(hc)-abs(hd));
end
```

**Figure 4.** Error function used within the Kanamori et al. (1999) method.

# STATE-SPACE APPROACH

The governing second-order differential equation of the SDOF system (Eq. (2)) can be written as a first-order continuous-time state-space model as follows:

$$\dot{x}(t) = \mathbf{A}x(t) + \mathbf{B}\ddot{u}_g(t), \qquad (41)$$

where $x(t) = [u(t) \quad \dot{u}(t)]^T$ is the state vector, and

$$\mathbf{A} = \begin{bmatrix} 0 & 1 \\ -\omega_n^2 & -2\xi\omega_n \end{bmatrix}, \qquad (42)$$

$$\mathbf{B} = \begin{bmatrix} 0 \\ -1 \end{bmatrix}. \qquad (43)$$

The state vector might be a hidden variable and not directly measurable (see, e.g., Williams and Lawrence, 2007), so there is a second equation describing the measurement/observation, which is expressed as

$$u(t) = \mathbf{C}x(t) + \mathbf{D}\ddot{u}_g(t). \qquad (44)$$

For the relative displacement calculation in this paper, $\mathbf{C} = [1 \quad 0]$ and $\mathbf{D} = 0$.

To use this state-space model for dynamic analysis under arbitrary ground motion excitation, it needs to be converted to a discrete-time version as:

$$x[k+1] = \mathbf{A}_D x[k] + \mathbf{B}_D \ddot{u}_g[k] \qquad (45)$$

$$u[k] = \mathbf{C}_D x[k] + \mathbf{D}_D \ddot{u}_g[k], \qquad (46)$$

through which the relative displacement can be easily calculated recursively, similar to Eq. (10). By solving Eq. (41), we can write

$$x((k+1)\Delta t) = e^{\mathbf{A}\Delta t}x(k\Delta t) + e^{\mathbf{A}(k+1)\Delta t} \mathbf{B} \int_{k\Delta t}^{(k+1)\Delta t} e^{-\mathbf{A}\tau}\ddot{u}_g(t)d\tau, \qquad (47)$$

from which the exact discrete-time equivalence of matrix $\mathbf{A}$ is $\mathbf{A}_D = e^{\mathbf{A}\Delta t}$. Depending on how the ground motion acceleration is assumed to vary in each time interval, two different methods emerge as discussed earlier in the transfer function section, through which two different sets of discrete-time versions are calculated for the remaining matrices.

## *Hold Methods*

One of the most common approaches in solving state-space equations is the ZOH assumption, in which it is assumed that the input excitation, here ground motion acceleration, remains constant during each time interval, i.e., $\ddot{u}_g(t) = \ddot{u}_g[k]$ for $k\Delta t \leq t < (k+1)\Delta t$. Using this assumption, we have:

$$\mathbf{B}_D = e^{\mathbf{A}(k+1)\Delta t}\mathbf{B} \int_{k\Delta t}^{(k+1)\Delta t} e^{-\mathbf{A}\tau} d\tau = \mathbf{A}^{-1}(\mathbf{A}_D - \mathbf{I})\mathbf{B}, \tag{48}$$

and trivially,

$$\mathbf{C}_D = \mathbf{C} \tag{49}$$

$$\mathbf{D}_D = \mathbf{D}. \tag{50}$$

Therefore, if the input excitation is truly stair-like, this discrete-time state-space model is an exact representation of the continuous-time state-space model.

Similar to the transfer function section, we can assume that the ground acceleration varies linearly between consecutive samples, that is,

$$\ddot{u}_g(t) = \ddot{u}_g(k\Delta t) + \frac{t-k\Delta t}{\Delta t}\left[\ddot{u}_g((k+1)\Delta t) - \ddot{u}_g(k\Delta t)\right] \quad \text{for} \quad k\Delta t \leq t < (k+1)\Delta t. \tag{51}$$

Replacing Eq. (51) into Eq. (47), the discrete-time matrices are obtained as follows (see Franklin et al., 1997, and Hoagg, 2004 for detailed derivation)

$$\mathbf{A}_D = e^{\mathbf{A}\Delta t} \tag{52}$$

$$\mathbf{B}_D = \frac{1}{\Delta t}\mathbf{A}^{-2}(\mathbf{A}_D - \mathbf{I})^2 \mathbf{B} \tag{53}$$

$$\mathbf{C}_D = \mathbf{C} \tag{54}$$

$$\mathbf{D}_D = \mathbf{D} + \mathbf{C}\left[\frac{1}{\Delta t}\mathbf{A}^{-2}(\mathbf{A}_D - \mathbf{I}) - \mathbf{A}^{-1}\right]\mathbf{B}. \tag{55}$$

*Approximate Methods*

The matrix exponential is a computationally expensive operation. To simplify the calculation, numerical integration methods can be used to approximate the solution of the integral rather than using the exact expression. Such approximations, and their effects on other matrices, lead to various discretization methods. The derivation of these matrices becomes more straightforward if we first convert the continuous-time state-space equations to the Laplace domain and then use the correspondence between the s- and z-domains, as will be discussed below. Finally, by applying the inverse z-transform, we can obtain the discrete-time state-space equations.

By applying the Laplace transform to the continuous-time state-space equations, we obtain:

$$s\mathbf{X}(s) = \mathbf{A}\mathbf{X}(s) + \mathbf{B}\ddot{U}_g(s), \tag{56}$$

$$\mathbf{U}(s) = \mathbf{C}\mathbf{X}(s) + \mathbf{D}\ddot{U}_g(s). \tag{57}$$

In the following subsections, we use the correspondence between s- and z-domains, which was presented in the transfer function section, and we obtain corresponding discrete-time state space equations.

**Euler Methods**

To convert the continuous-time state space model to the discrete-time model using the Forward Euler method, we replace $s$ with $\frac{z-1}{\Delta t}$ in Eqs. (56) and (57), and we use the z-Transform properties. It is straightforward to show that

$$\mathbf{A}_D = \mathbf{I} + \mathbf{A}\Delta t \tag{58}$$

$$\mathbf{B}_D = \mathbf{B}\Delta t \tag{59}$$

$$\mathbf{C}_D = \mathbf{C} \tag{60}$$

$$\mathbf{D}_D = \mathbf{D}. \tag{61}$$

In the Backward Euler method, we use $s \approx \frac{z-1}{z\Delta t}$, which converts the state Eq. (56) to

$$\frac{z-1}{z\Delta t}X(s) = \mathbf{A}X(s) + \mathbf{B}\ddot{U}_g(s), \tag{62}$$

which corresponds to the time-domain equation

$$x[k+1] - x[k] = \mathbf{A}\Delta t x[k+1] + \mathbf{B}\Delta t \ddot{u}_g[k+1]. \tag{63}$$

By rearranging Eq. (63), we have

$$x[k+1] - \mathbf{A}\Delta t x[k+1] - \mathbf{B}\Delta t \ddot{u}_g[k+1] = x[k]. \tag{64}$$

We define $v[k+1] = x[k]$, so we have

$$(\mathbf{I} - \mathbf{A}\Delta t)x[k] = v[k] + \mathbf{B}\Delta t \ddot{u}_g[k], \tag{65}$$

which results in

$$x[k] = (\mathbf{I} - \mathbf{A}\Delta t)^{-1}v[k] + (\mathbf{I} - \mathbf{A}\Delta t)^{-1}\mathbf{B}\Delta t \ddot{u}_g[k]. \tag{66}$$

By replacing $x[k]$ with $v[k+1]$, we can obtain the standard form of the discrete-time state space model as

$$v[k+1] = (\mathbf{I} - \mathbf{A}\Delta t)^{-1}v[k] + (\mathbf{I} - \mathbf{A}\Delta t)^{-1}\mathbf{B}\Delta t \ddot{u}_g[k]. \tag{67}$$

Therefore,

$$\mathbf{A}_D = (\mathbf{I} - \mathbf{A}\Delta t)^{-1} \tag{68}$$

$$\mathbf{B}_D = (\mathbf{I} - \mathbf{A}\Delta t)^{-1}\mathbf{B}\Delta t. \tag{69}$$

By replacing Eq. (66) in the observation equation (Eq. (44)), we will have

$$u[k] = \mathbf{C}\{(\mathbf{I} - \mathbf{A}\Delta t)^{-1}v[k] + (\mathbf{I} - \mathbf{A}\Delta t)^{-1}\mathbf{B}\Delta t \ddot{u}_g[k]\} + \mathbf{D}\ddot{u}_g[k], \tag{70}$$

therefore

$$\mathbf{C}_D = \mathbf{C}(\mathbf{I} - \mathbf{A}\Delta t)^{-1} \tag{71}$$

$$\mathbf{D}_D = \mathbf{D} + \mathbf{C}(\mathbf{I} - \mathbf{A}\Delta t)^{-1}\mathbf{B}\Delta t. \tag{72}$$

**Tustin Method**

To discretize the state-space model using the Tustin method, we replace $s$ in Eq. (56) with $\frac{2(z-1)}{\Delta t(z+1)}$. After converting to the time-domain, we have

$$x[k+1] - x[k] = \frac{A\Delta t}{2}(x[k+1] + x[k]) + \frac{B\Delta t}{2}(\ddot{u}_g[k+1] + \ddot{u}_g[k]). \tag{73}$$

Following the same approach we used for the backward Euler method, we move all $k+1$ terms to the left-hand side, so

$$\left(I - \frac{A\Delta t}{2}\right)x[k+1] - \frac{B\Delta t}{2}\ddot{u}_g[k+1] = \left(I + \frac{A\Delta t}{2}\right)x[k] + \frac{B\Delta t}{2}\ddot{u}_g[k]. \tag{74}$$

By defining $v[k+1] = \left(I + \frac{A\Delta t}{2}\right)x[k] + \frac{B\Delta t}{2}\ddot{u}_g[k]$, we have

$$\left(I - \frac{A\Delta t}{2}\right)x[k] = v[k] + \frac{B\Delta t}{2}\ddot{u}_g[k], \tag{75}$$

from which

$$x[k] = \left(I - \frac{A\Delta t}{2}\right)^{-1} v[k] + \left(I - \frac{A\Delta t}{2}\right)^{-1} \frac{B\Delta t}{2}\ddot{u}_g[k]. \tag{76}$$

Now, by replacing Eq. (76) into the definition of $v[k+1]$, we have

$$v[k+1] = \left(I + \frac{A\Delta t}{2}\right)\left(I - \frac{A\Delta t}{2}\right)^{-1} v[k] + \left(I - \frac{A\Delta t}{2}\right)^{-1} B\Delta t \ddot{u}_g[k]. \tag{77}$$

Therefore,

$$A_D = \left(I + \frac{A\Delta t}{2}\right)\left(I - \frac{A\Delta t}{2}\right)^{-1} \tag{78}$$

$$B_D = \left(I - \frac{A\Delta t}{2}\right)^{-1} B\Delta t. \tag{79}$$

To obtain the observation equation, we substitute Eq. (76) into Eq. (44), and we have

$$u[k] = C\left\{\left(I - \frac{A\Delta t}{2}\right)^{-1} v[k] + \left(I - \frac{A\Delta t}{2}\right)^{-1} \frac{B\Delta t}{2}\ddot{u}_g[k]\right\} + D\ddot{u}_g[k], \tag{80}$$

therefore

$$C_D = C\left(I - \frac{A\Delta t}{2}\right)^{-1} \tag{81}$$

$$D_D = D + C\left(I - \frac{A\Delta t}{2}\right)^{-1} \frac{B\Delta t}{2}. \tag{82}$$

For readers' convenience, I have summarized all discrete-time matrices obtained from the various discretization methods in Table 11. Additionally, a MATLAB function implementing all these methods for time-history analysis is presented in **Figure 5**.

Table 11. Discrete-time state-space matrices obtained using various discretization methods.

|       | ZOH | FOH | FE | BE | Tustin |
|-------|-----|-----|----|----|--------|
| $A_D$ | $e^{A\Delta t}$ | $e^{A\Delta t}$ | $I + A\Delta t$ | $(I - A\Delta t)^{-1}$ | $\left(I + \frac{A\Delta t}{2}\right)\left(I - \frac{A\Delta t}{2}\right)^{-1}$ |
| $B_D$ | $A^{-1}(A_D - I)B$ | $\frac{1}{\Delta t}A^{-2}(A_D - I)^2 B$ | $B\Delta t$ | $(I - A\Delta t)^{-1}B\Delta t$ | $\left(I - \frac{A\Delta t}{2}\right)^{-1} B\Delta t$ |
| $C_D$ | C | C | C | $C(I - A\Delta t)^{-1}$ | $C\left(I - \frac{A\Delta t}{2}\right)^{-1}$ |
| $D_D$ | D | $D + C\left[\frac{1}{\Delta t}A^{-2}(A_D - I) - A^{-1}\right]B$ | D | $D + C(I - A\Delta t)^{-1}B\Delta t$ | $D + C\left(I - \frac{A\Delta t}{2}\right)^{-1}\frac{B\Delta t}{2}$ |

```matlab
function [disp,Time] = StateSpace(wn,xi,ag,dt,dta,dtd,c2d_method,interp_flag)
% Inputs:
%           % wn: natural frequency (rad/sec)
%           % xi: damping ratio
%           % ag: ground acceleration (Nx1)
%           % dt: ground motion dt (sec.)
%           % dta: analysis dt (sec.) if interpolation flag is set to 1 or 2
%           % dtd: display dt (sec.)
%           % c2d_method: method used for time discretization
%           % interp_flag: interpolation flag (0: no interpolation, 1: linear, 2: sinc)
% Outputs:
%           % disp : relative displacement response (Nx1)
%           % Time : time vector (Nx1)

if interp_flag==1 % linear interpolation
        US = dt/dta;
        N = length(ag);
        Time = [0:dt:(N-1)*dt]';
        ta = [0:dta:(US*N-1)*dta]';
        ag = interp1(Time, ag, ta, 'linear', 'extrap');
        dt = dta;
elseif interp_flag==2
        US = dt/dta;
        ag = interp(ag,US);
        dt = dta;
end
N = length(ag);
Time = [0:dt:(N-1)*dt]';

A = [0,1;-wn^2,-2*xi*wn];
B = [0;-1];
C = [1,0];
D = 0;
wd = wn*sqrt(1-xi^2);
if strcmp(c2d_method,'ZOH')
        Ad = exp(-xi*wn*dt)*[cos(wd*dt)+xi/sqrt(1-xi^2)*sin(wd*dt),1/wd*sin(wd*dt);-wn/sqrt(1-xi^2)*sin(wd*dt),cos(wd*dt)-xi/sqrt(1-xi^2)*sin(wd*dt)];
        Bd = A\(Ad-eye(2))*B;
        Cd = C;
        Dd = D;
elseif strcmp(c2d_method,'FOH')
        Ad = exp(-xi*wn*dt)*[cos(wd*dt)+xi/sqrt(1-xi^2)*sin(wd*dt),1/wd*sin(wd*dt);-wn/sqrt(1-xi^2)*sin(wd*dt),cos(wd*dt)-xi/sqrt(1-xi^2)*sin(wd*dt)];
        Bd = 1/dt*A^(-2)*(Ad-eye(2))^2*B;
        Cd = C;
        Dd = D+C*(1/dt*A^(-2)*(Ad-eye(2))-inv(A))*B;
elseif strcmp(c2d_method,'FE')
        Ad = eye(2)+A*dt;
        Bd = B*dt;
        Cd = C;
        Dd = D;
elseif strcmp(c2d_method,'BE')
        Ad = eye(2)/(eye(2)-A*dt);
        Bd = Ad*B*dt;
        Cd = C*Ad;
        Dd = D+Cd*B*dt;
elseif strcmp(c2d_method,'Tustin')
        Ad = (eye(2)+0.5*A*dt)/(eye(2)-0.5*A*dt);
        Bd = (eye(2)-0.5*A*dt)\B*dt;
        Cd = C/(eye(2)-0.5*A*dt);
        Dd = D+C/(eye(2)-0.5*A*dt)*B*dt/2;
end

disp = zeros(N,1);
x = [0;0];
for k=2:N
        x = Ad*x+Bd*ag(k-1);
        disp(k,1) = Cd*x+Dd*ag(k);
end

if dt>dtd
        US = dt/dtd;
        disp = interp(disp,US);
        Time = [0:dtd:(US*N-1)*dtd]';
end
end
```

**Figure 5.** State-space MATLAB function for calculating the time-history response of an SDOF system.

## EARTHQUAKE ENGINEERING METHODS

In the previous sections, two classical families of methods for analyzing LTI systems, widely used for many years across various fields, were reviewed, and the corresponding formula for analyzing an SDOF system was presented. In earthquake engineering, the numerical solution to Eq. (2) has been calculated through time-stepping methods since the seminal work by Newmark (1959). These methods can be generally classified into three groups (Chopra, 2007b): 1- methods based on interpolation of the excitation, 2- methods based on finite difference expression of velocity and acceleration, and 3- methods based on assumed variation of acceleration. As will be discussed further, while these methods might have been developed independently from classical solutions presented earlier, they can be classified into one of those groups, except the Newmark method (Newmark, 1959), which does not directly fit into those discretization methods.

## *Methods Based on Excitation Interpolation*

As the analytical solution to step and linear excitations is possible, one approach is to assume input excitation varies linearly between two samples and solve Eq. (2) under initial conditions at the beginning of each time step, step excitation due to the excitation at the beginning of the time step, and linear excitation due to the variation within the time step. By adding these three responses, the total response of the system at the end of each time step is calculated and is used as the initial conditions for the next time step. Nigam and Jennings (1968, 1969) were among the first to propose this method. According to this method, relative displacement and velocity at a discrete time instant $k+1$ are calculated using displacement and velocity from the previous step and ground acceleration at the current and next time steps (Nigam and Jennings, 1969)

$$\begin{bmatrix} u[k+1] \\ \dot{u}[k+1] \end{bmatrix} = \widehat{\mathbf{A}} \begin{bmatrix} u[k] \\ \dot{u}[k] \end{bmatrix} + \widehat{\mathbf{B}} \begin{bmatrix} \ddot{u}_g[k] \\ \ddot{u}_g[k+1] \end{bmatrix}. \tag{83}$$

with $\widehat{\mathbf{A}} = \begin{bmatrix} a_{11} & a_{12} \\ a_{21} & a_{22} \end{bmatrix}$ and $\widehat{\mathbf{B}} = \begin{bmatrix} b_{11} & b_{12} \\ b_{21} & b_{22} \end{bmatrix}$ where

$$a_{11} = e^{-\xi \omega_n \Delta t} \left[ \frac{\xi}{\sqrt{1-\xi^2}} \sin(\omega_d \Delta t) + \cos(\omega_d \Delta t) \right] \tag{84}$$

$$a_{12} = \frac{e^{-\xi \omega_n \Delta t}}{\omega_d} \sin(\omega_d \Delta t) \tag{85}$$

$$a_{21} = -\frac{\omega_n}{\sqrt{1-\xi^2}} e^{-\xi \omega_n \Delta t} \sin(\omega_d \Delta t) \tag{86}$$

$$a_{22} = e^{-\xi \omega_n \Delta t} \left[ \cos(\omega_d \Delta t) - \frac{\xi}{\sqrt{1-\xi^2}} \sin(\omega_d \Delta t) \right] \tag{87}$$

$$b_{11} = e^{-\xi \omega_n \Delta t} \left[ \left( \frac{2\xi^2-1}{\omega_n^2 \Delta t} + \frac{\xi}{\omega_n} \right) \frac{\sin(\omega_d \Delta t)}{\omega_d} + \left( \frac{2\xi}{\omega_n^3 \Delta t} + \frac{1}{\omega_n^2} \right) \cos(\omega_d \Delta t) \right] - \frac{2\xi}{\omega_n^3 \Delta t} \tag{88}$$

$$b_{12} = -e^{-\xi \omega_n \Delta t} \left[ \left( \frac{2\xi^2-1}{\omega_n^2 \Delta t} \right) \frac{\sin(\omega_d \Delta t)}{\omega_d} + \frac{2\xi}{\omega_n^3 \Delta t} \cos(\omega_d \Delta t) \right] - \frac{1}{\omega_n^2} + \frac{2\xi}{\omega_n^3 \Delta t} \tag{89}$$

$$b_{21} = e^{-\xi \omega_n \Delta t} \left\{ \left( \frac{2\xi^2-1}{\omega_n^2 \Delta t} + \frac{\xi}{\omega_n} \right) \left[ \cos(\omega_d \Delta t) - \frac{\xi}{\sqrt{1-\xi^2}} \sin(\omega_d \Delta t) \right] - \left( \frac{2\xi}{\omega_n^3 \Delta t} + \frac{1}{\omega_n^2} \right) [\omega_d \sin(\omega_d \Delta t) + \xi \omega_n \cos(\omega_d \Delta t)] \right\} + \frac{1}{\omega_n^2 \Delta t} \tag{90}$$

$$b_{22} = -e^{-\xi \omega_n \Delta t} \left\{ \frac{2\xi^2-1}{\omega_n^2 \Delta t} \left[ \cos(\omega_d \Delta t) - \frac{\xi}{\sqrt{1-\xi^2}} \sin(\omega_d \Delta t) \right] - \frac{2\xi}{\omega_n^3 \Delta t} [\omega_d \sin(\omega_d \Delta t) + \xi \omega_n \cos(\omega_d \Delta t)] \right\} - \frac{1}{\omega_n^2 \Delta t}. \tag{91}$$

This method has been used by strong motion networks for several years (see, e.g., Wald et al., 2021). However, looking at Eq. (83), this is a discrete-time state-space model, and based on the assumption of the input excitation variation, which is a linear variation, it seems to be the FOH discretization approach already discussed in previous section. To verify this fact, let's rewrite the standard state-space model of Eqs. (45) and (46), and add relative velocity too. That is,

$$\boldsymbol{x}[k+1] = \mathbf{A}_D \boldsymbol{x}[k] + \mathbf{B}_D \ddot{u}_g[k] \tag{92}$$

$$\begin{bmatrix} u[k] \\ \dot{u}[k] \end{bmatrix} = \mathbf{C}_D \boldsymbol{x}[k] + \mathbf{D}_D \ddot{u}_g[k], \tag{93}$$

in which $\mathbf{A}_D$, $\mathbf{B}_D$, $\mathbf{C}_D$, and $\mathbf{D}_D$ are calculated according to Table 11, where **A** and **B** are still the same, but $\mathbf{C} = \mathbf{I}_{2 \times 2}$, and $\mathbf{D} = [0 \quad 0]^T$. It is important to note that while the state vector $\boldsymbol{x}(t) = [u(t) \quad \dot{u}(t)]^T$ in the continuous-time domain has physical meaning, once the system is converted to the discrete-time domain through discretization, $x[k]$ may no longer have the same physical meaning. So, Eq. (92) cannot be compared to Eq. (83). However, the output of the observation equation remains unchanged. Writing the output equation for time step $k + 1$, we have

$$\begin{bmatrix} u[k+1] \\ \dot{u}[k+1] \end{bmatrix} = \mathbf{C}_D x[k+1] + \mathbf{D}_D \ddot{u}_g[k+1]. \tag{94}$$

Now, replacing $x[k + 1]$ from Eq. (92), we have

$$\begin{bmatrix} u[k+1] \\ \dot{u}[k+1] \end{bmatrix} = \mathbf{C}_D (\mathbf{A}_D x[k] + \mathbf{B}_D \ddot{u}_g[k]) + \mathbf{D}_D \ddot{u}_g[k+1]. \tag{95}$$

Using Eq. (93), we can write

$$x[k] = \mathbf{C}_D^{-1} \begin{bmatrix} u[k] \\ \dot{u}[k] \end{bmatrix} - \mathbf{C}_D^{-1} \mathbf{D}_D \ddot{u}_g[k]. \tag{96}$$

Now, by substituting $x[k]$ into Eq. (95), we have

$$\begin{bmatrix} u[k+1] \\ \dot{u}[k+1] \end{bmatrix} = \mathbf{C}_D \left\{ \mathbf{A}_D \left( \mathbf{C}_D^{-1} \begin{bmatrix} u[k] \\ \dot{u}[k] \end{bmatrix} - \mathbf{C}_D^{-1} \mathbf{D}_D \ddot{u}_g[k] \right) + \mathbf{B}_D \ddot{u}_g[k] \right\} + \mathbf{D}_D \ddot{u}_g[k+1], \tag{97}$$

which can be simplified as follows using the fact that $\mathbf{C}_D = \mathbf{C} = \mathbf{I}$ in the FOH discretization method

$$\begin{bmatrix} u[k+1] \\ \dot{u}[k+1] \end{bmatrix} = \mathbf{A}_D \begin{bmatrix} u[k] \\ \dot{u}[k] \end{bmatrix} + [\mathbf{B}_D - \mathbf{A}_D \mathbf{D}_D \quad \mathbf{D}_D] \begin{bmatrix} \ddot{u}_g[k] \\ \ddot{u}_g[k+1] \end{bmatrix}. \tag{98}$$

Eq. (98) is identical to Eq. (83) (Nigam and Jennings, 1969) provided that $\mathbf{A}_D = \widehat{\mathbf{A}}$ and $[\mathbf{B}_D - \mathbf{A}_D \mathbf{D}_D \quad \mathbf{D}_D] = \widehat{\mathbf{B}}$. The continuous-time matrix $\mathbf{A} = \begin{bmatrix} 0 & 1 \\ -\omega_n^2 & -2\xi\omega_n \end{bmatrix}$ has complex-conjugate eigenvalues as $-\xi\omega_n \mp i\omega_d$, so according to the Cayley-Hamilton theorem (see, e.g., Walden and Roelof, 1982), we have

$$\mathbf{A}_D = e^{\mathbf{A}\Delta t} = e^{-\xi\omega_n \Delta t} \left[ \cos(\omega_d \Delta t) \mathbf{I} + \frac{1}{\omega_d} \sin(\omega_d \Delta t) (\mathbf{A} + \xi\omega_n \mathbf{I}) \right], \tag{99}$$

which can be expanded as

$$e^{\mathbf{A}\Delta t} = e^{-\xi\omega_n \Delta t} \begin{bmatrix} \cos(\omega_d \Delta t) + \frac{\xi}{\sqrt{1-\xi^2}} \sin(\omega_d \Delta t) & \frac{1}{\omega_d} \sin(\omega_d \Delta t) \\ -\frac{\omega_n}{\sqrt{1-\xi^2}} \sin(\omega_d \Delta t) & \cos(\omega_d \Delta t) - \frac{\xi}{\sqrt{1-\xi^2}} \sin(\omega_d \Delta t) \end{bmatrix}, \tag{100}$$

which is identical to matrix $\widehat{\mathbf{A}}$. It is also straightforward, but cumbersome, to show that $[\mathbf{B}_D - \mathbf{A}_D \mathbf{D}_D \quad \mathbf{D}_D] = \widehat{\mathbf{B}}$. Therefore, the method proposed by Nigam and Jennings is identical to the discrete-time state-space model discretized using the FOH method.

A MATLAB function for calculating the time-history response of an SDOF system using the method of Nigam and Jennings (1969) is provided in **Figure 6**, although an equivalent response can be obtained using the state-space method with FOH discretization (**Figure 5**).

```matlab
function [disp,Time] = Nigam_Jennings1969(wn,xi,ag,dt,dta,dtd,interp_flag)
% Inputs:
        % wn: natural frequency (rad/sec)
        % xi: damping ratio
        % ag: ground acceleration (Nx1)
        % dt: ground motion dt (sec.)
        % dta: analysis dt (sec.) if interpolation flag is set to 1 or 2
        % dtd: display dt (sec.)
        % interp_flag: interpolation flag (0: no interpolation, 1: linear, 2: sinc)

% Outputs:
        % disp : relative displacement response (Nx1)
        % Time : time vector (Nx1)

if interp_flag==1
        US = dt/dta;
        N = length(ag);
        Time = [0:dt:(N-1)*dt]';
        ta = [0:dta:(US*N-1)*dta]';
        ag = interp1(Time,ag,ta,'linear', 'extrap');
        dt = dta;
elseif interp_flag==2
        US = dt/dta;
        ag = interp(ag,US);
        dt = dta;
end
N = length(ag);
Time = [0:dt:(N-1)*dt]';

wd = wn*sqrt(1-xi^2);
alpha = exp(-xi*wn*dt);
theta = wd*dt;
b = xi/sqrt(1-xi^2);
a11 = alpha*(b*sin(theta)+cos(theta));
a12 = alpha/wd*sin(theta);
a21 = -wn/sqrt(1-xi^2)*alpha*sin(theta);
a22 = alpha*(cos(theta)-b*sin(theta));
f1 = (2*xi^2-1)/(wn^2*dt)+xi/wn;
f2 = 2*xi/(wn^3*dt)+1/wn^2;
b11 = alpha*(f1*sin(theta)/wd+f2*cos(theta))-2*xi/(wn^3*dt);
f3 = (2*xi^2-1)/(wn^2*dt);
f4 = 2*xi/(wn^3*dt);
b12 = -alpha*(f3*sin(theta)/wd+f4*cos(theta))-1/wn^2+2*xi/(wn^3*dt);
b21 = alpha*(f1*(cos(theta)-b*sin(theta))-f2*(wd*sin(theta)+xi*wn*cos(theta)))+1/(wn^2*dt);
b22 = -alpha*(f3*(cos(theta)-b*sin(theta))-f4*(wd*sin(theta)+xi*wn*cos(theta)))-1/(wn^2*dt);

A = [a11,a12;a21,a22];
B = [b11,b12;b21,b22];

state = zeros(2,N);
for k=1:N-1
        state(:,k+1) = A*state(:,k)+B*[ag(k);ag(k+1)];
end

disp = state(1,:)';
if dt>dtd
        US = dt/dtd;
        disp = interp(disp,US);
        Time = [0:dtd:(US*N-1)*dtd]';
end
```

**Figure 6.** A MATLAB function for calculating the time-history response of an SDOF system using the method of Nigam and Jennings (1969).

## *Methods Based on Finite Difference of Velocity and Acceleration*

This method is commonly used in most finite element software for dynamic analysis (Bathe, 1982). In this approach, velocity and acceleration are approximated as

$$\dot{u}[k] = \frac{u[k+1]-u[k-1]}{2\Delta t} \tag{101}$$

$$\ddot{u}[k] = \frac{u[k+1]-2u[k]+u[k-1]}{\Delta t^2}. \tag{102}$$

By substituting these expressions into Eq. (2), we obtain

$$u[k+1] = -\frac{\Delta t^2}{1+\xi\omega_n\Delta t}\ddot{u}_g[k] - \frac{\omega_n^2\Delta t^2 - 2}{1+\xi\omega_n\Delta t}u[k] - \frac{1-\xi\omega_n\Delta t}{1+\xi\omega_n\Delta t}u[k-1], \qquad (103)$$

Through this relation, the displacement at the next time step is calculated from the displacements at the current and previous time steps. Comparing Eq. (103) with Eq. (10), this method can be classified as a transfer function method, with coefficients reported in Table 12.

A MATLAB function for calculating the time-history response of an SDOF system using the method Central Difference is provided in **Figure 7**, although an equivalent response can be obtained using the transfer function method with CD discretization (**Figure 2**).

Table 12. Coefficients of the discrete-time transfer function discretized using the Central Difference method.

| Coefficient | Value |
| --- | --- |
| $a_1$ | $(\omega_n^2\Delta t^2 - 2)/(1+\xi\omega_n\Delta t)$ |
| $a_2$ | $(1-\xi\omega_n\Delta t)/(1+\xi\omega_n\Delta t)$ |
| $b_0$ | 0 |
| $b_1$ | $-\Delta t^2/(1+\xi\omega_n\Delta t)$ |
| $b_2$ | 0 |

```matlab
function [disp,Time] = CentralDiff(wn,xi,ag,dt,dta,dtd,interp_flag)
% Inputs:
%         % wn: natural frequency (rad/sec)
%         % xi: damping ratio
%         % ag: ground acceleration (Nx1)
%         % dt: ground motion dt (sec.)
%         % dta: analysis dt (sec.) if interpolation flag is set to 1 or 2
%         % dtd: display dt (sec.)
%         % interp_flag: interpolation flag (0: no interpolation, 1: linear, 2: sinc)

% Outputs:
%         % disp : relative displacement response (Nx1)
%         % Time : time vector (Nx1)

if interp_flag==1
        US = dt/dta;
        N = length(ag);
        Time = [0:dt:(N-1)*dt]';
        ta = [0:dta:(US*N-1)*dta]';
        ag = interp1(Time,ag,ta,'linear', 'extrap');
        dt = dta;
elseif interp_flag==2
        US = dt/dta;
        ag = interp(ag,US);
        dt = dta;
end
N = length(ag);
Time = [0:dt:(N-1)*dt]';

khat = 1/dt^2+2*xi*wn/(2*dt);
a = 1/dt^2-2*xi*wn/(2*dt);
b = wn^2-2/dt^2;
disp = zeros(N,1);
vel = zeros(N,1);
acc = zeros(N,1);
acc(1,1) = -ag(1)-(2*xi*wn)*vel(1,1)-wn^2*disp(1,1);
u0 = disp(1,1)-dt*vel(1,1)+dt^2/2*acc(1,1);
for k=1:N-1
        if k==1
                p_hat = -ag(k)-a*u0-b*disp(k,1);
                disp(k+1,1) = p_hat/khat;
                vel(k,1) = (disp(k+1,1)-u0)/(2*dt);
                acc(k,1) = (disp(k+1)-2*disp(k,1)+u0)/dt^2;
        else
                p_hat = -ag(k)-a*disp(k-1,1)-b*disp(k,1);
                disp(k+1,1) = p_hat/khat;
                vel(k,1) = (disp(k+1,1)-disp(k-1,1))/(2*dt);
                acc(k,1) = (disp(k+1)-2*disp(k,1)+disp(k-1,1))/dt^2;
        end
end

if dt>dtd
US = dt/dtd;
disp = interp(disp,US);
Time = [0:dtd:(US*N-1)*dtd]';
end
```

**Figure 7.** A MATLAB function for calculating the time-history response of an SDOF system using the Central Difference method.

## Methods Based on Assumed Acceleration Variation

While methods presented earlier can be classified under classical method, the most common technique in earthquake engineering, proposed by Newmark (1959), does not fall into either of these classes. A brief summary of this method is provided below.

The velocity of a system at the end of a time interval can be written as

$$\dot{u}[k+1] = \dot{u}[k] + \int_{t[k]}^{t[k]+\Delta t} \ddot{u}(\tau) \, d\tau. \tag{104}$$

To evaluate the above integral, the functional form of the acceleration within the interval $\Delta t$ must be known. Instead, Newmark approximated the integral (area under the acceleration curve) by a rectangle with height equal to the weighted average of the accelerations at the two ends of the interval. Accordingly,

$$\dot{u}[k+1] = \dot{u}[k] + (1-\gamma)\Delta t \ddot{u}[k] + \gamma \Delta t \ddot{u}[k+1]. \tag{105}$$

A common assumption is $\gamma = 0.5$, which assigns equal weight to both ends. To obtain displacement, consider a Taylor expansion:

$$u[k+1] = u[k] + \Delta t \dot{u}[k] + \frac{1}{2}\Delta t^2 \ddot{u}[k] + \cdots. \tag{106}$$

If acceleration were constant, the first three terms would suffice. However, since acceleration varies during each interval, Newmark retained the first three terms but replaced the coefficient of $\frac{1}{2}\ddot{u}[k]$ with a weighted average similar to the velocity approximation:

$$u[k+1] = u[k] + \Delta t \dot{u}[k] + \Delta t^2 \left[\left(\frac{1}{2} - \beta\right)\Delta t \ddot{u}[k] + \beta \ddot{u}[k+1]\right], \tag{107}$$

where $\frac{1}{6} \leq \beta \leq \frac{1}{4}$. The upper limit $\beta = \frac{1}{4}$ gives equal weight to both ends and is referred to as the *constant average acceleration method*, while the lower limit $\beta = \frac{1}{6}$ assigns twice the weight to $\ddot{u}[k]$ compared to $\ddot{u}[k+1]$, leading to the *linear acceleration method*.

From Eq. (105), the velocity increment is

$$\Delta \dot{u}[k] = \dot{u}[k+1] - \dot{u}[k] = \Delta t \ddot{u}[k] + \gamma \Delta t \Delta \ddot{u}[k], \tag{108}$$

where $\Delta \ddot{u}[k] = \ddot{u}[k+1] - \ddot{u}[k]$. Similarly, from Eq. (107), the displacement increment is

$$\Delta u[k] = u[k+1] - u[k] = \Delta t \dot{u}[k] + \frac{\Delta t^2}{2}\ddot{u}[k] + \beta \Delta t^2 \Delta \ddot{u}[k], \tag{109}$$

from which it follows that

$$\Delta \ddot{u}[k] = \frac{1}{\beta \Delta t^2}\Delta u[k] - \frac{1}{\beta \Delta t}\dot{u}[k] - \frac{1}{2\beta}\ddot{u}[k]. \tag{110}$$

Substituting Eq. (110) into Eq. (108) yields

$$\Delta \dot{u}[k] = \frac{\gamma}{\beta \Delta t}\Delta u[k] - \frac{\gamma}{\beta}\dot{u}[k] + \Delta t \left(1 - \frac{\gamma}{2\beta}\right)\ddot{u}[k]. \tag{111}$$

Since the equation of motion (Eq. (2)) holds at each step,

$$\Delta \ddot{u}[k] + 2\xi \omega_n \Delta \dot{u}[k] + \omega_n^2 \Delta u[k] = -\left(\ddot{u}_g[k+1] - \ddot{u}_g[k]\right). \tag{112}$$

Replacing $\Delta \ddot{u}[k]$ and $\Delta \dot{u}[k]$ from Eqs. (110) and (111), the displacement increment can be expressed as

$$\Delta u[k] = \frac{\Delta \hat{P}}{\hat{K}}. \tag{113}$$

where

$$\Delta \hat{P} = -\ddot{u}_g[k+1] + \ddot{u}_g[k] - a\dot{u}[k] - b\ddot{u}[k] \tag{114}$$

$$\hat{K} = \omega_n^2 + \frac{\gamma}{\beta \Delta t}2\xi \omega_n + \frac{1}{\beta \Delta t^2}, \tag{115}$$

with

$$a = \frac{1}{\beta \Delta t} + \frac{\gamma}{\beta}2\xi \omega_n \tag{116}$$

$$b = \frac{1}{2\beta} + \Delta t \left(\frac{\gamma}{2\beta} - 1\right)2\xi \omega_n. \tag{117}$$

Having obtained $\Delta u[k]$, the increments $\Delta \ddot{u}[k]$ and $\Delta \dot{u}[k]$ are calculated from Eqs. (110) and (111). Consequently, the displacement, velocity, and acceleration at the next step are updated as (Chopra, 2007b):

$$u[k+1] = u[k] + \Delta u[k] \tag{118}$$

$$\dot{u}[k+1] = \dot{u}[k] + \Delta \dot{u}[k] \tag{119}$$

$$\ddot{u}[k+1] = \ddot{u}[k] + \Delta \ddot{u}[k]. \tag{120}$$

A MATLAB function for calculating the time-history response of an SDOF system using the Newmark method is provided in **Figure 8**. In the Newmark method, acceleration is treated as an additional state, so the method cannot be expressed in a standard state-space form. However, by defining the state vector as $\boldsymbol{x}_n = [u[k], \dot{u}[k], \ddot{u}[k]]^T$, Eqs. (105), (107), and (2) can be written in a matrix–vector form as

$$\mathbf{H}_1 \boldsymbol{x}[k+1] = \mathbf{H}_0 \boldsymbol{x}[k] + \mathbf{F} \ddot{u}_g[k+1], \tag{121}$$

where

$$\mathbf{H}_1 = \begin{bmatrix} 1 & 0 & -\beta \Delta t^2 \\ 0 & 1 & -\gamma \Delta t \\ \omega_n^2 & 2\xi \omega_n & 1 \end{bmatrix} \tag{122}$$

$$\mathbf{H}_0 = \begin{bmatrix} 1 & \Delta t & \left(\frac{1}{2} - \beta\right) \Delta t^2 \\ 0 & 1 & (1-\gamma)\Delta t \\ 0 & 0 & 0 \end{bmatrix} \tag{123}$$

$$\mathbf{F} = [0 \quad 0 \quad -1]^T. \tag{124}$$

Eq. (121) can be converted into a standard discrete-time state-space model for displacement calculation as

$$\boldsymbol{x}[k+1] = \mathbf{A}_D \boldsymbol{x}[k] + \mathbf{B}_D \ddot{u}_g[k+1].^{\dagger} \tag{125}$$

$$u[k] = \mathbf{C}_D \boldsymbol{x}[k] + \mathbf{D}_D \ddot{u}_g[k]. \tag{126}$$

where $\mathbf{C}_D = [1 \quad 0 \quad 0]$, $\mathbf{D}_D = 0$, and

$$\mathbf{A}_D = \frac{1}{1 + 2\xi \omega_n \gamma \Delta t + \omega_n^2 \beta \Delta t^2} \begin{bmatrix} a_{11} & a_{12} & a_{13} \\ a_{21} & a_{22} & a_{23} \\ a_{31} & a_{32} & a_{33} \end{bmatrix} \tag{127}$$

$$\mathbf{B}_D = \frac{1}{1 + 2\xi \omega_n \gamma \Delta t + \omega_n^2 \beta \Delta t^2} [-\beta \Delta t^2 \quad -\gamma \Delta t \quad -1]^T \tag{128}$$

in which

$$a_{11} = 1 + 2\xi \omega_n \gamma \Delta t \tag{129}$$

$$a_{12} = -2\xi \omega_n \beta \Delta t^2 + 2\xi \omega_n \gamma \Delta t^2 + \Delta t \tag{130}$$

$$a_{13} = -2\xi \omega_n \beta \Delta t^3 + \xi \omega_n \gamma \Delta t^3 - \beta \Delta t^2 + \frac{1}{2} \Delta t^2 \tag{131}$$

$$a_{21} = -\omega_n^2 \gamma \Delta t \tag{132}$$

---

[†] Note that the input excitation appears with a forward time delay.

$$a_{22} = 1 + \omega_n^2 \beta \Delta t^2 - \omega_n^2 \gamma \Delta t^2 \tag{133}$$

$$a_{23} = \omega_n^2 \beta \Delta t^3 - \frac{1}{2} \omega_n^2 \gamma \Delta t^3 - \gamma \Delta t + \Delta t \tag{134}$$

$$a_{31} = -\omega_n^2 \tag{135}$$

$$a_{32} = -2\xi\omega_n - \omega_n^2 \Delta t \tag{136}$$

$$a_{33} = \omega_n^2 \beta \Delta t^2 - \frac{1}{2} \omega_n^2 \Delta t^2 + 2\xi\omega_n\gamma\Delta t - 2\xi\omega_n\Delta t. \tag{137}$$

Therefore, the Newmark method can also be interpreted as a state-space model. This alternative representation is particularly useful for investigating the method's stability later in this paper.

```matlab
function [disp,Time] = Newmark(wn,xi,ag,dt,dta,dtd,gamma,beta,interp_flag)
% Inputs:
        % wn: natural frequency (rad/sec)
        % xi: damping ratio
        % ag: ground acceleration (Nx1)
        % dt: ground motion dt (sec.)
        % dta: analysis dt (sec.) if interpolation flag is set to 1 or 2
        % dtd: display dt (sec.)
        % gamma: Newmark's integration parameter
        % beta: Newmark's integration parameter
        % interp_flag: interpolation flag (0: no interpolation, 1: linear, 2: sinc)
% Outputs:
        % disp : relative displacement response (Nx1)
        % Time : time vector (Nx1)

if interp_flag==1
        US = dt/dta;
        N = length(ag);
        Time = [0:dt:(N-1)*dt]';
        ta = [0:dta:(US*N-1)*dta]';
        ag = interp1(Time,ag,ta,'linear', 'extrap');
        dt = dta;
elseif interp_flag==2
        US = dt/dta;
        ag = interp(ag,US);
        dt = dta;
end
N = length(ag);
Time = [0:dt:(N-1)*dt]';

khat = wn^2+gamma/(beta*dt)*(2*xi*wn)+1/(beta*dt^2);
a = 1/(beta*dt)+gamma/beta*(2*xi*wn);
b = 1/(2*beta)+dt*(gamma/(2*beta)-1)*(2*xi*wn);
disp = zeros(N,1);
vel = zeros(N,1);
acc = zeros(N,1);
acc(1,1) = -ag(1)-(2*xi*wn)*vel(1)-wn^2*disp(1);
for k=1:N-1
        deltap_hat = -ag(k+1)-(-ag(k))+a*vel(k)+b*acc(k);
        deltau = deltap_hat/khat;
        deltav = gamma/(beta*dt)*deltau-gamma/beta*vel(k)+dt*(1-gamma/(2*beta))*acc(k);
        deltaa = 1/(beta*dt^2)*deltau-1/(beta*dt)*vel(k)-1/(2*beta)*acc(k);
        disp(k+1) = disp(k)+deltau;
        vel(k+1) = vel(k)+deltav;
        acc(k+1) = acc(k)+deltaa;
end

if dt>dtd
        US = dt/dtd;
        disp = interp(disp,US);
        Time = [0:dtd:(US*N-1)*dtd]';
end
end
```

**Figure 8.** A MATLAB function for calculating the time-history response of an SDOF system using the Newmark method.

# STABILITY STUDY

*Transfer function Models*

A system is stable if any input signal with bounded amplitude applied to the system produces an output that also remains bounded for all time. This type of stability is called Bounded-Input, Bounded-Output (BIBO) stability (Oppenheim et al., 1997). A continuous-time transfer function is BIBO stable if all its poles lie strictly in the left half of the complex s-plane (Oppenheim et al., 1997). The SDOF system with the transfer function of Eq. (6) is strictly stable as long as it has positive damping because its poles are $p_{1,2} = -\xi\omega_n \mp i\omega_d$. However, when converted to a discrete-time transfer function, the stability may or may not be preserved. A discrete-time transfer function is BIBO stable if all its poles are inside the unit circle (Oppenheim, 1999). All discrete-time transfer functions introduced in this paper are in the form of

$$H(z) = \frac{b_0 + b_1 z^{-1} + b_2 z^{-2}}{1 + a_1 z^{-1} + a_2 z^{-2}}. \tag{138}$$

The poles of this generic transfer function are calculated by finding the roots of the denominator

$$p_{1,2} = \frac{-a_1 \mp \sqrt{a_1^2 - 4a_2}}{2}. \tag{139}$$

I have calculated the poles of transfer functions discretized through all methods presented in this paper, which are reported in Table 13. As seen in this table, except for the Forward Euler, Least-Squares (LSQ), and the method proposed by Kanamori et al. (1999), all other discretization methods preserve the stability of the continuous-time transfer function. As discussed earlier, the LSQ method, if solved using the iterative nonlinear optimization algorithm, enforces the stability of the system during the process, so the final system would be stable.

Table 13. Stability conditions of Transfer Function methods.

| Discretization Method | Poles | \|Poles\| | Stability |
|---|---|---|---|
| ZOH | $e^{-\xi\omega_n \Delta t}[\cos(\omega_d \Delta t) \mp i \sin(\omega_d \Delta t)]$ | $e^{-\xi\omega_n \Delta t}$ | Stable |
| FOH | $e^{-\xi\omega_n \Delta t}[\cos(\omega_d \Delta t) \mp i \sin(\omega_d \Delta t)]$ | $e^{-\xi\omega_n \Delta t}$ | Stable |
| Impulse | $e^{-\xi\omega_n \Delta t}[\cos(\omega_d \Delta t) \mp i \sin(\omega_d \Delta t)]$ | $e^{-\xi\omega_n \Delta t}$ | Stable |
| FE | $1 - \xi\omega_n \Delta t \mp i\omega_d \Delta t$ | $\sqrt{1 - 2\xi\omega_n \Delta t + \omega_n^2 \Delta t^2}$ | $\Delta t < \frac{2\xi}{\omega_n}$ |
| BE | $\frac{1 + \xi\omega_n \Delta t \mp i\omega_d \Delta t}{1 + 2\xi\omega_n \Delta t + \omega_n^2 \Delta t^2}$ | $1/\sqrt{1 + 2\xi\omega_n \Delta t + \omega_n^2 \Delta t^2}$ | Stable |
| Tustin | $\frac{4 - \omega_n^2 \Delta t^2 \mp 4i\omega_d \Delta t}{4 + 4\xi\omega_n \Delta t + \omega_n^2 \Delta t^2}$ | $\frac{\sqrt{4 - 4\xi\omega_n \Delta t + \omega_n^2 \Delta t^2}}{\sqrt{4 + 4\xi\omega_n \Delta t + \omega_n^2 \Delta t^2}}$ | Stable |
| Matched | $e^{-\xi\omega_n \Delta t \mp i\omega_d \Delta t}$ | $e^{-\xi\omega_n \Delta t}$ | Stable |
| LSQ* | $\left(-a_1 \mp \sqrt{a_1^2 - 4a_2}\right)/2$ | $\left\|-a_1 \mp \sqrt{a_1^2 - 4a_2}\right\|/2$ | $\left\|-a_1 \mp \sqrt{a_1^2 - 4a_2}\right\| < 2$ |
| Kanamori et al. | $\frac{1 + \xi\omega_n \Delta t \mp i\omega_d \Delta t}{1 + 2\xi\omega_n \Delta t + \omega_n^2 \Delta t^2}$ | $1/\sqrt{1 + 2\xi\omega_n \Delta t + \omega_n^2 \Delta t^2}$ | $-\frac{\omega_n \Delta t}{2} < \xi < 1$ |

* The stability of this method without the iterative optimization solution is not guaranteed.

*State-Space Models*

To study the necessary conditions for the stability of the state-space models, we first convert them to the transfer function representation. Using Eq. (56), we can obtain the Laplace Transform of the state vector as $X(s) = (s\mathbf{I} - \mathbf{A})^{-1}\mathbf{B}\ddot{U}_g(s)$. By replacing this relationship in Eq. (57), we have

$$\frac{U(s)}{\ddot{U}_g(s)} = [\mathbf{C}(s\mathbf{I} - \mathbf{A})^{-1}\mathbf{B} + \mathbf{D}], \tag{140}$$

which is the continuous-time relative displacement and velocity Transfer Function of the system. Using the definition of the matrix inversion, we can rewrite Eq. (140) as

$$\frac{U(s)}{\ddot{U}_g(s)} = \frac{\mathbf{C}\,adj(s\mathbf{I}-\mathbf{A})\,\mathbf{B}+\mathbf{D}|s\mathbf{I}-\mathbf{A}|}{|s\mathbf{I}-\mathbf{A}|}, \tag{141}$$

where $adj(.)$ denotes the adjoint of a matrix. According to the definition of stability of a transfer function discussed earlier, the state-space system is stable if all poles of this transfer function lie strictly in the left half of the complex s-plane. The roots of the $|s\mathbf{I} - \mathbf{A}| = 0$ are eigenvalues of the matrix $\mathbf{A}$. So, the stability condition implies that all eigenvalues must have a negative real part. It is easy to show that the eigenvalues of matrix $\mathbf{A}$ are $s_{1,2} = -\xi\omega_n \mp i\omega_d$, which are identical to the poles of the continuous-time transfer function and have a negative real part (when damping is positive), so the state-space system in continuous-time is always stable. For the discrete-time state-space model (Eqs. (45) and (46)), it is easy to show that

$$\frac{U(z)}{\ddot{U}_g(z)} = [\mathbf{C}_D(z\mathbf{I} - \mathbf{A}_D)^{-1}\mathbf{B}_D + \mathbf{D}_D], \tag{142}$$

which can be rewritten as

$$\frac{U(z)}{\ddot{U}_g(z)} = \frac{\mathbf{C}_D\,adj(z\mathbf{I}-\mathbf{A}_D)\,\mathbf{B}_D+\mathbf{D}_D|z\mathbf{I}-\mathbf{A}_D|}{|z\mathbf{I}-\mathbf{A}_D|}. \tag{143}$$

Now, according to the stability condition of a discrete-time transfer function, the discrete-time state-space model is stable if the eigenvalues of the matrix $\mathbf{A}_D$ are all inside the unit circle. I have calculated the eigenvalues of matrix $\mathbf{A}_D$ for all discretization methods, and collected them in Table 14. Not surprisingly, eigenvalues in this table are identical to the poles presented in Table 13, and so are the stability conditions.

Table 14. Stability conditions of state-space methods.

| Discretization Method | Eigenvalues | \|Eigenvalues\| | Stability |
|---|---|---|---|
| ZOH | $e^{-\xi\omega_n\Delta t}[\cos(\omega_d\Delta t) \mp i\sin(\omega_d\Delta t)]$ | $e^{-\xi\omega_n\Delta t}$ | Stable |
| FOH | $e^{-\xi\omega_n\Delta t}[\cos(\omega_d\Delta t) \mp i\sin(\omega_d\Delta t)]$ | $e^{-\xi\omega_n\Delta t}$ | Stable |
| FE | $1 - \xi\omega_n\Delta t \mp i\omega_d\Delta t$ | $\sqrt{1 - 2\xi\omega_n\Delta t + \omega_n^2\Delta t^2}$ | $\Delta t < \frac{2\xi}{\omega_n}$ |
| BE | $\frac{1+\xi\omega_n\Delta t \mp i\omega_d\Delta t}{1+2\xi\omega_n\Delta t+\omega_n^2\Delta t^2}$ | $1/\sqrt{1 + 2\xi\omega_n\Delta t + \omega_n^2\Delta t^2}$ | Stable |
| Tustin | $\frac{4-\omega_n^2\Delta t^2 \mp 4i\omega_d\Delta t}{4+4\xi\omega_n\Delta t+\omega_n^2\Delta t^2}$ | $\frac{\sqrt{4-4\xi\omega_n\Delta t+\omega_n^2\Delta t^2}}{\sqrt{4+4\xi\omega_n\Delta t+\omega_n^2\Delta t^2}}$ | Stable |

*Earthquake Engineering Methods*

The stability condition of the method proposed by Nigam and Jennings (1968, 1969) is identical to that of the state-space model with FOH discretization; therefore, the method is unconditionally stable, as shown in Table 14. To study the stability of the Central Difference method, we can take advantage of its transfer function representation (Table 12). The poles of the discrete-time transfer function are

$$p_{1,2} = \frac{2 - \omega_n^2 \Delta t^2 \mp i\omega_n \Delta t \sqrt{4(1-\xi^2) - \omega_n^2 \Delta t^2}}{2(1+\xi\omega_n \Delta t)}. \tag{144}$$

These poles lie inside the unit circle unconditionally because

$$|p_{1,2}| = \frac{\sqrt{1 - \xi^2 \omega_n^2 \Delta t^2}}{1 + \xi\omega_n \Delta t} \leq 1. \tag{145}$$

Therefore, the Central Difference method is unconditionally stable provided that the poles are complex-conjugate, as shown in Eq. (144). If the poles are real, they are given by

$$p_{1,2} = \frac{2 - \omega_n^2 \Delta t^2 \mp \omega_n \Delta t \sqrt{\omega_n^2 \Delta t^2 - 4(1-\xi^2)}}{2(1+\xi\omega_n \Delta t)}. \tag{146}$$

For the general case of $\xi = 0$, the absolute value of the largest pole is

$$|p_{1,2}| = \frac{2 - \omega_n^2 \Delta t^2 + \omega_n \Delta t \sqrt{\omega_n^2 \Delta t^2 - 4}}{2}, \tag{147}$$

which lies inside the unit circle when $\Delta t < 2/\omega_n$.

To investigate the stability of the Newmark method, we take advantage of its state-space representation. Specifically, the method is stable if the eigenvalues of matrix $\mathbf{A}_D$, defined in Eq. (127), all lie inside the unit circle, which is equivalent to having a spectral radius less than 1, as extensively discussed by Nickell (1971) and Bathe and Wilson (1973). For $\gamma = 0.5$, the nontrivial eigenvalues of matrix $\mathbf{A}_D$ are shown below:

$$\lambda_{1,2} = \frac{1 + b_1\left(\beta - \frac{1}{2}\right) \pm i\sqrt{b_1^2\left(\beta - \frac{1}{4}\right) + b_1 - \frac{b_2^2}{4}}}{1 + \frac{b_2}{2} + \beta b_1}. \tag{148}$$

where the two dimensionless parameters are $b_1 = \Delta t^2 \omega_n^2$ and $b_2 = 2\xi\Delta t\omega_n$. The absolute value of these Eigenvalues are

$$|\lambda_{1,2}| = \sqrt{\frac{1 + \frac{b_2}{2} - \beta b_1}{1 + \frac{b_2}{2} + \beta b_1}}, \tag{149}$$

which are always less than 1. However, similar to the Central Difference method, Eq. (149) is valid if $b_1^2\left(\beta - \frac{1}{4}\right) + b_1 - \frac{b_2^2}{4} \geq 0$ which gives $\Delta t < \frac{2\sqrt{\frac{1-\xi^2}{1-4\beta}}}{\omega_n}$. For the critical case of $\xi = 0$, the stability condition reduces to $\Delta t < \frac{2}{\omega_n\sqrt{1-4\beta}}$ which is the expression presented in textbooks (e.g., Chopra, 2007). Therefore, the Newmark method with the constant average acceleration assumption ($\beta = 1/4$) is unconditionally stable, while the linear acceleration approach is stable if $\Delta t/T_n \leq 0.551$.

## IMPLEMENTATION REMARKS

In real-life applications, ground motions are typically recorded at sampling rates of 100 to 200 Hz. However, the accuracy—and in some cases, even the stability—of the numerical methods discussed in the previous sections is sensitive to the choice of $\Delta t$ relative to the natural period of the SDOF, $T_n = 2\pi/\omega_n$. A smaller analysis step size, $\Delta t$, makes the discrete-time system more closely approximate its continuous-time counterpart. Therefore, when performing time-history analysis, it is recommended to use a reduced step size, $\Delta t_a$, instead of the original sampling interval, where the upsampling factor is defined as $\Delta t / \Delta t_a \geq 1$.

Once the continuous-time system is discretized in a smaller sampling time, $\Delta t_a$, the ground motion has to be upsampled as well to be able to use it within the recursive calculation. Linear interpolation is the most common approach, especially in earthquake engineering, so it will be employed in this paper. However, linear interpolation is not a perfect reconstruction method of a signal and might miss some between-sample peaks (details are beyond the scope of this paper). So, in addition to linear interpolation, we run examples using Sinc interpolation of ground motions, which is a perfect way of reconstructing a continuous-time signal (ground truth signal) using discrete-time samples (see, e.g., Oppenheim, 1999).

Finally, after the system response is computed—whether at the original ground motion sampling interval, $\Delta t$, or at the analysis step size, $\Delta t_a$—it is often useful to interpolate the response. This is particularly important if the results are to be compared with a ground-truth analytical solution, which is typically represented on a computer using a very small sampling interval, $\Delta T$, or if a time-domain metric, such as the peak response, is required. The interpolation should be performed using the Sinc method to ensure that between-sample peaks are accurately captured. Accordingly, I set $\Delta t_d = \Delta T$, where $\Delta t_d$ is the sampling interval used for displaying the calculated responses.

## EXAMPLES

The purpose of this section is not to recommend which method should be used. The performance of each method would depend on the objective (e.g., peak response is of interest or the entire time-history), on the application (e.g., if a delay in calculation or response is fine), and on the factors like stability or complexity. Herein, I show the performance of all discussed methods for a few examples for which analytical solutions are available, and to see how the accuracy can change depending on factors like natural frequency of the SDOF, sampling time, and interpolation technique.

As the first example, consider a ground motion acceleration as a pure sinusoidal function $\ddot{u}_g(t) = \ddot{u}_g^o \sin \omega_0 t$ where $\ddot{u}_g^o = 1$ is the amplitude and $\omega_0$ is the frequency of excitation in rad/sec. The analytical solution is

$$u(t) = e^{-\xi \omega_n \Delta t}(C_4 \sin \omega_d t + C_3 \cos \omega_d t) + C_1 \sin \omega_0 t + C_2 \cos \omega_0 t, \tag{150}$$

where

$$C_1 = -\frac{\ddot{u}_g^o}{\omega_n^2} \frac{1 - \left(\frac{\omega_0}{\omega_n}\right)^2}{\left[1 - \left(\frac{\omega_0}{\omega_n}\right)^2\right]^2 + \left[2\xi\left(\frac{\omega_0}{\omega_n}\right)\right]^2} \tag{151}$$

$$C_2 = -\frac{\ddot{u}_g^o}{\omega_n^2} \frac{-2\xi\left(\frac{\omega_0}{\omega_n}\right)}{\left[1 - \left(\frac{\omega_0}{\omega_n}\right)^2\right]^2 + \left[2\xi\left(\frac{\omega_0}{\omega_n}\right)\right]^2} \tag{152}$$

$$C_3 = \frac{\ddot{u}_g^o}{\omega_n^2} \frac{-2\xi\left(\frac{\omega_0}{\omega_n}\right)}{\left[1 - \left(\frac{\omega_0}{\omega_n}\right)^2\right]^2 + \left[2\xi\left(\frac{\omega_0}{\omega_n}\right)\right]^2} \tag{153}$$

$$C_4 = \frac{\xi \omega_n C_3 - C_1 \omega_0}{\omega_d}. \tag{154}$$

To display this continuous-time signal, $\Delta T = 0.0001\ sec.$ is used. Also, as mentioned earlier, after calculating the response through numerical methods, the response is Sinc-interpolated using the same sampling time used for the continuous-time signal, i.e., $\Delta t_d = \Delta T$, to make sure two signals are compatible in time.

**Figure 9** shows the performance of three time-stepping methods commonly used in the earthquake engineering community for a system with a natural period of $T_n = 0.3\ sec.$ and a damping ratio of $\xi = 5\%$, subjected to harmonic excitation with a frequency five times the natural frequency. The discrete-time input excitation is sampled at 100 Hz, and the same sampling interval is used in the analysis, i.e., $\Delta t_a = \Delta t = 0.01\ sec$. As shown in the figure, small errors appear in the calculated responses, among which the solution obtained using the Central Difference method is slightly more accurate. Note that $\beta = 1/6$ is used for the Newmark method.

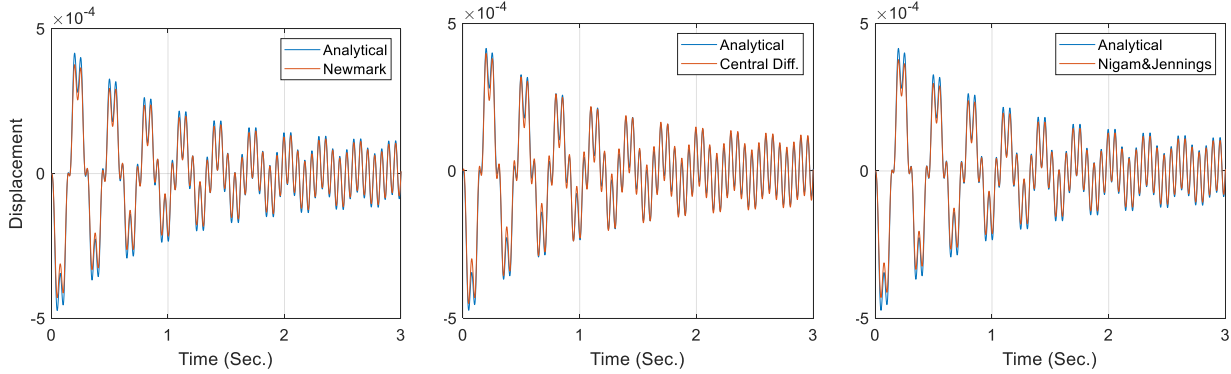

**Figure 9.** Comparison between Earthquake Engineering methods and the analytical solution ($T_n = 0.3\ sec., \xi = 5\%, \omega_0 = 5\omega_n, \Delta t = 0.01\ sec.$, no interpolation)

**Figure 10** shows a similar comparison for transfer function methods. The Forward Euler method is excluded because it is unstable according to Table 13. Among the remaining methods, the ZOH, Impulse, and Matched methods exhibit better performance, although a small phase shift is observed in the ZOH and Matched methods due to input excitation discretization and the presence of zeros at infinity, respectively. For this example, the Nyquist frequency is 50 Hz while the natural frequency is 3.33 Hz; therefore, only limited improvement is expected from the Tustin method with pre-warping compared to the case without pre-warping. The difference between these two methods becomes more evident in **Figure 11**, where the pre-warped Tustin method is able to accurately match the dominant (natural) frequency of the frequency response function. The final observation from **Figure 10** concerns Kanamori's method. Since the input excitation has a frequency around 17 Hz, a frequency bandwidth between 0 and 25 Hz is used for the minimization. As shown in **Figure 12**, while the amplitude of the frequency response function closely matches that of the ground-truth continuous-time system, the phase is not captured accurately because it is not included in the optimization. Consequently, the predicted response is not accurate, as clearly observed in **Figure 10**. The estimated effective natural periods and damping ratios are $\tilde{T}_n = 0.3$ and $\tilde{\xi} = -5.49\%$, respectively. In contrast, the LSQ method captures both amplitude and phase, resulting in a response that is almost as accurate as the other well-performing methods.

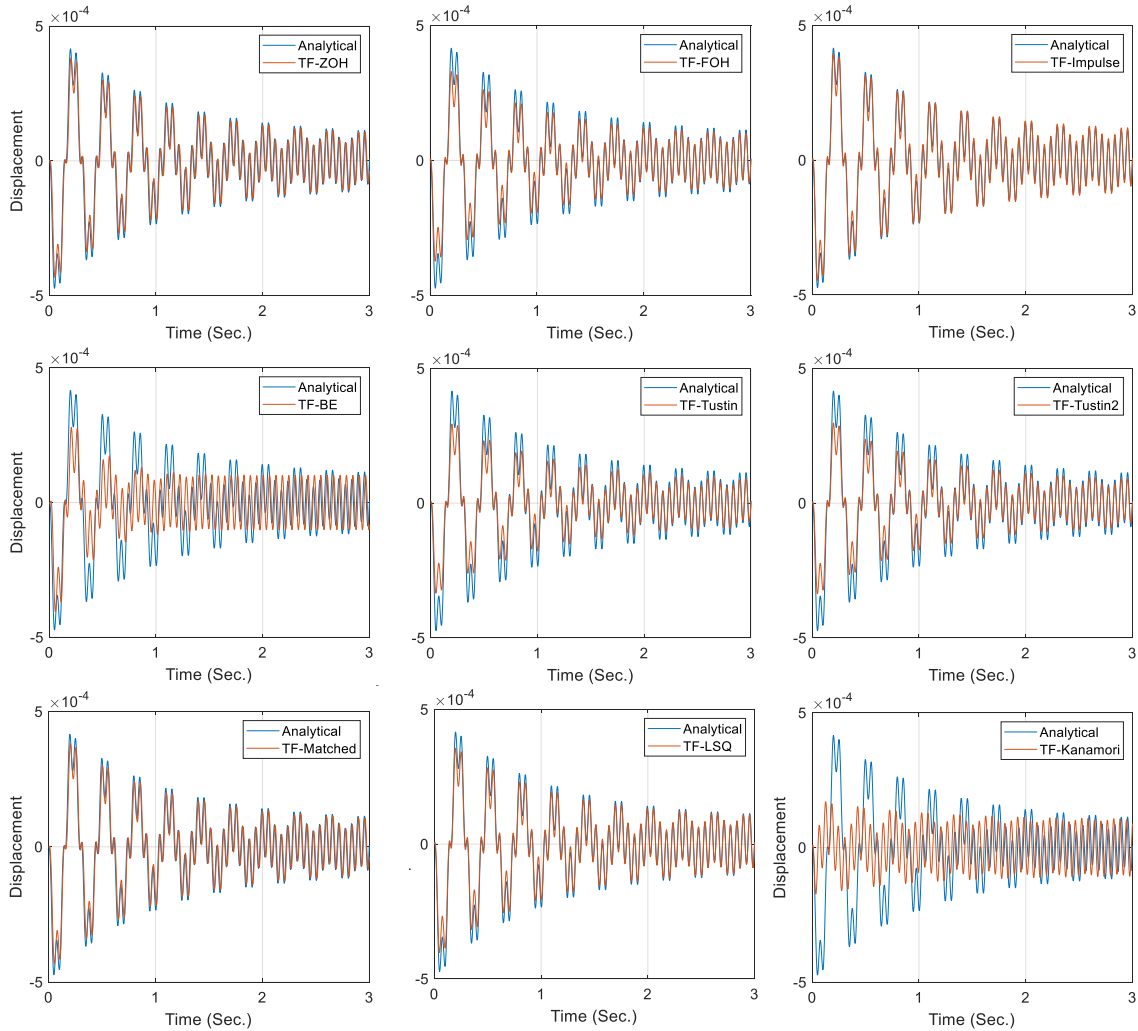

**Figure 10.** Comparison between transfer function methods and the analytical solution ($T_n = 0.3\ sec.$, $\xi = 5\%$, $\omega_0 = 5\omega_n$, $\Delta t = 0.01\ sec.$, no interpolation)

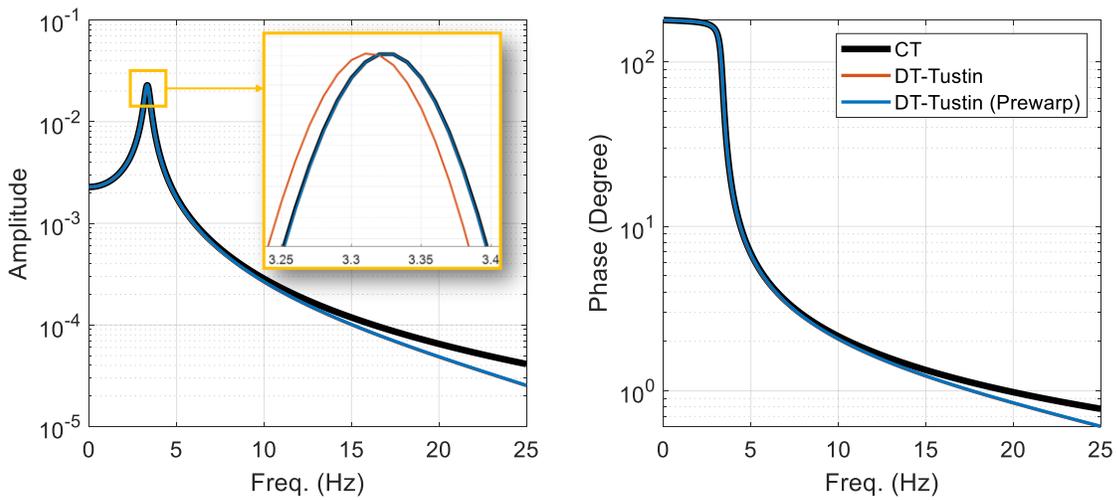

**Figure 11.** Comparison of amplitude (left) and phase (right) of the discrete-time frequency response functions estimated using the Tustin method, with and without pre-warping ($T_n = 0.3\ sec.$, $\xi = 5\%$, $\Delta t = 0.01\ sec.$, no interpolation)

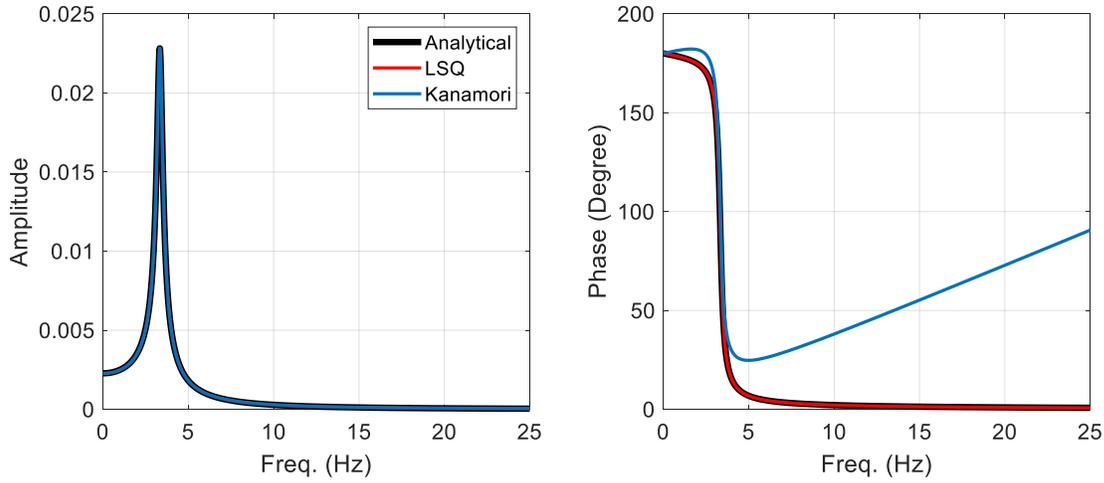

**Figure 12.** Comparison of amplitude (left) and phase (right) of the frequency response functions estimated using the LSQ and Kanamori et al. methods against the analytical method ($T_n = 0.3\ sec.$, $\xi = 5\%$, $\Delta t = 0.01\ sec.$, no interpolation)

The performance of the state-space methods is shown in **Figure 13**. The Forward Euler discretization method is omitted because it is unstable. Aside from the highly inaccurate results of the Backward Euler method, the performance of the other three methods is very similar, with the FOH method performing slightly better than the other two.

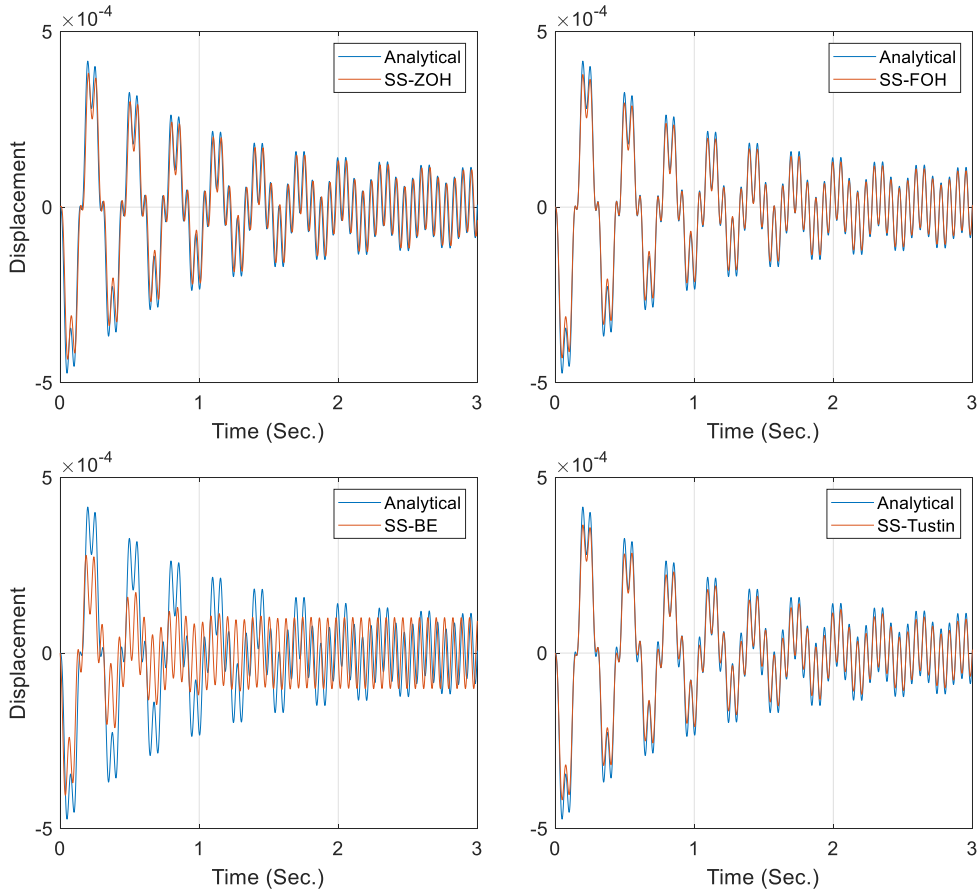

**Figure 13.** Comparison between state-space methods and the analytical solution ($T_n = 0.3\ sec.$, $\xi = 5\%$, $\omega_0 = 5\omega_n$, $\Delta t = 0.01\ sec.$, no interpolation)

To summarize the results of this example, **Figure 14** shows the Root-Mean-Square (RMS) errors and maximum displacement errors between all numerical methods and the analytical solution. The RMS and maximum displacement errors are defined as

$$\epsilon(S_d) = \frac{\max_k|u[k]| - \max_k|u_a[k]|}{\max_k|u_a[k]|} \times 100, \tag{155}$$

$$\epsilon(RMS) = \frac{\sqrt{\sum_k |u[k] - u_a[k]|^2}}{\sqrt{\sum_k |u_a[k]|^2}} \times 100, \tag{156}$$

where $u[k]$ and $u_a[k]$ are the numerically and analytically calculated relative displacements at the discrete-time instant $k\Delta t_d$. As seen in this figure, the performance of different methods may vary depending on the objective index. For example, the state-space model with FOH discretization performs among the best when RMS error is of interest, whereas there are other methods that perform better when peak displacement is considered. It is important to note, however, that the errors in all these methods are relatively large due to the analysis sampling time being large compared to the natural and excitation periods.

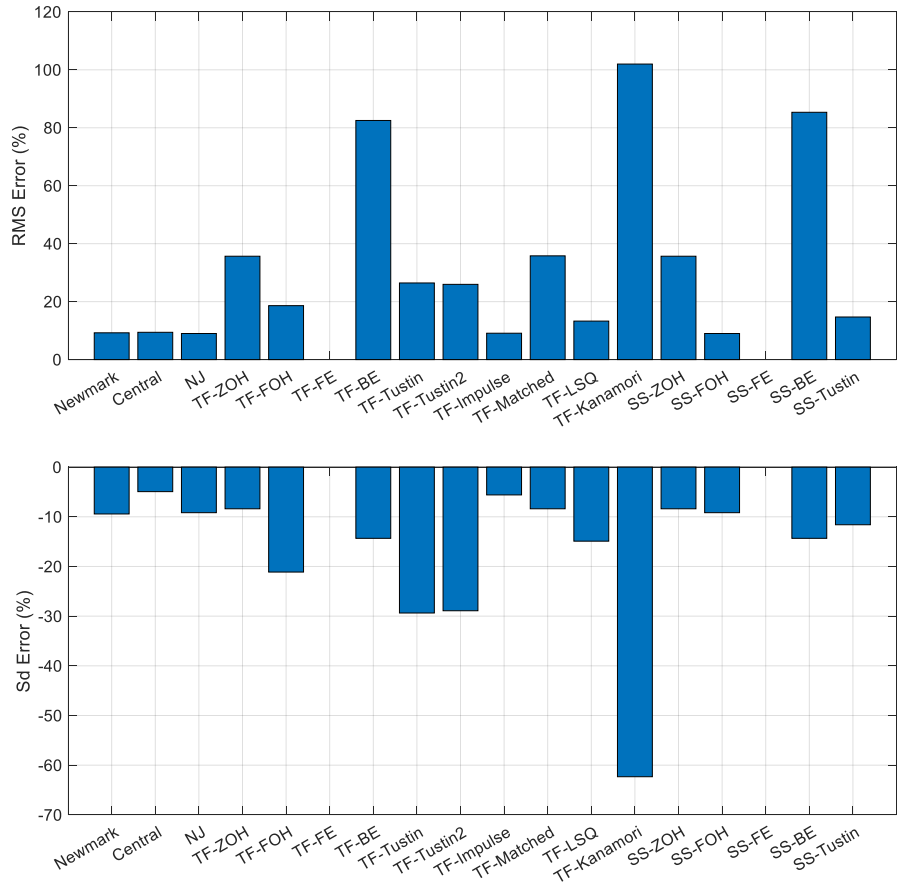

**Figure 14.** RMS and maximum displacement errors of various discretization methods ($T_n = 0.3\ sec., \xi = 5\%, \omega_0 = 5\omega_n, \Delta t = 0.01\ sec.$, no interpolation)

To improve the accuracy of the discretization methods, the same problem is solved using a sampling time 10 times smaller, i.e., $\Delta t_a = \frac{\Delta t}{10} = 0.001\ sec$. The input excitation, originally sampled at 100 Hz, is linearly interpolated to achieve this new sampling frequency. **Figure 15** shows the RMS and maximum displacement errors between all numerical methods and the analytical solution. With the smaller sampling time, the Forward Euler method becomes stable, so its results are also included in this figure. As seen, almost all methods now perform similarly, except for the Forward and Backward Euler methods, which still show larger errors. The performance of the Kanamori method is now comparable to that of the LSQ method because the discretization used to convert the continuous-time differential equation into a discrete-time difference equation is more accurate. This improvement can be seen more clearly in **Figure 16**. Although phase was not included in the optimization, the resulting phase is very close to the analytical and LSQ counterparts. Consequently, the calculated response is nearly identical to that obtained when phase is included in the optimization, as in the LSQ method. The estimated effective natural period and damping ratio are $\tilde{T}_n = 0.3$ and $\tilde{\xi} = 3.95\%$, respectively, which are very close to the actual properties of the SDOF system. It is important to emphasize that the level of errors (both RMS and maximum displacement) is still not acceptable for any of the methods. In other words, careful discretization of the system alone is not sufficient—the input excitation also needs to be accurately resampled.

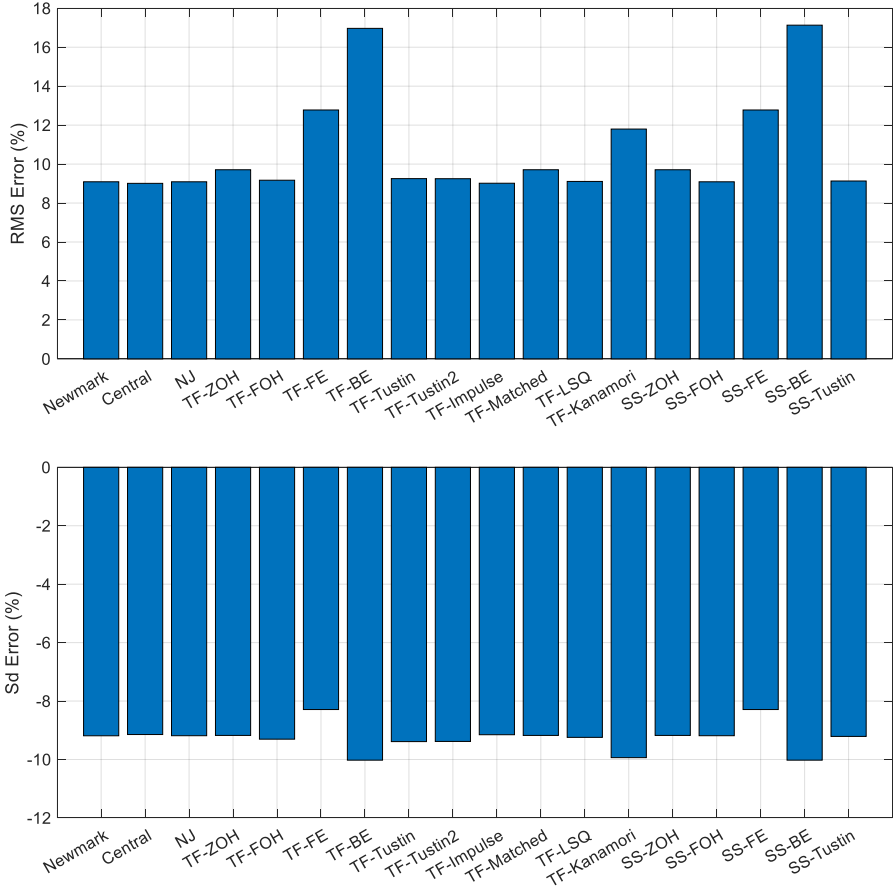

**Figure 15.** RMS and maximum displacement errors of various discretization methods ($T_n = 0.3\ sec.$, $\xi = 5\%$, $\omega_0 = 5\omega_n$, $\Delta t_a = 0.001\ sec.$, linear interpolation)

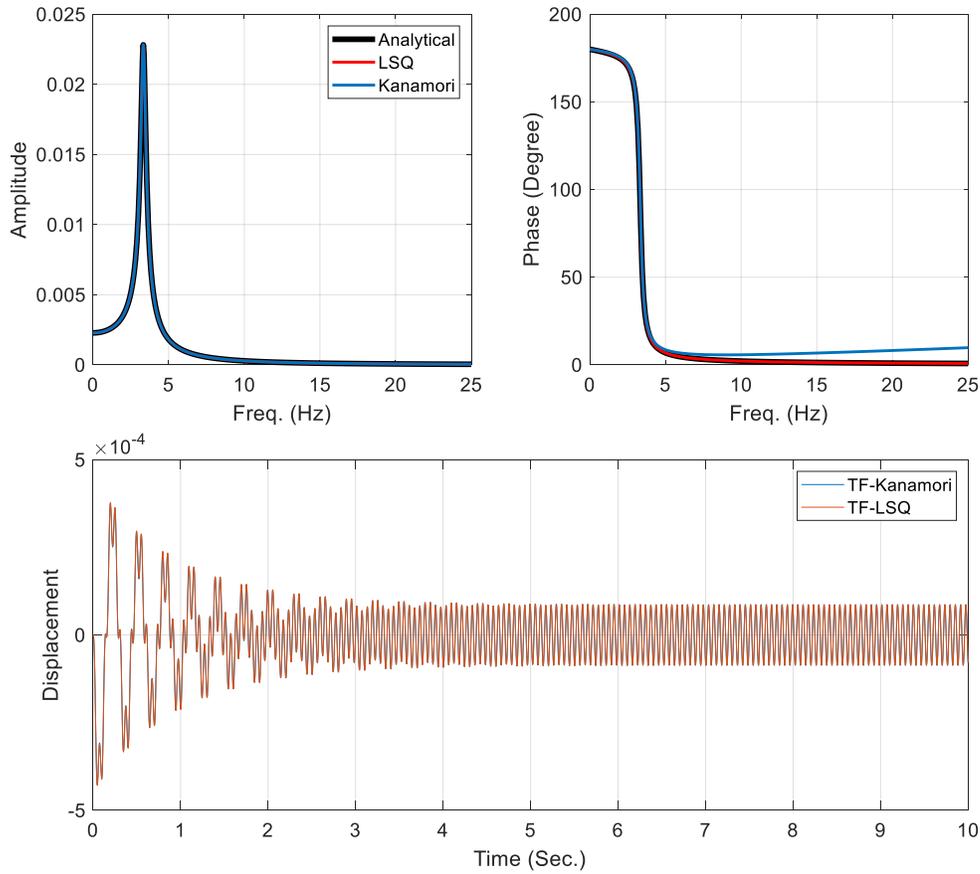

**Figure 16.** Top: Comparison of amplitude (left) and phase (right) of the frequency response functions estimated using the LSQ and Kanamori et al. methods against the analytical method. Bottom: Comparison between time-history responses obtained using LSQ and Kanamori methods ($T_n = 0.3\ sec., \xi = 5\%, \Delta t_a = 0.001\ sec.$, linear interpolation)

**Figure 17** shows the results obtained by changing the interpolation from linear to sinc interpolation. As seen, the error levels for almost all methods are very small. A few exceptions can be explained as follows. The change from linear to Sinc interpolation does not alter the discrete-time system properties and only affects the input excitation. Therefore, larger errors are still expected for the Euler methods, especially for the RMS index, where time synchronization is important. Time synchronization is also the source of the larger RMS errors observed in the ZOH, Matched, and Kanamori methods.

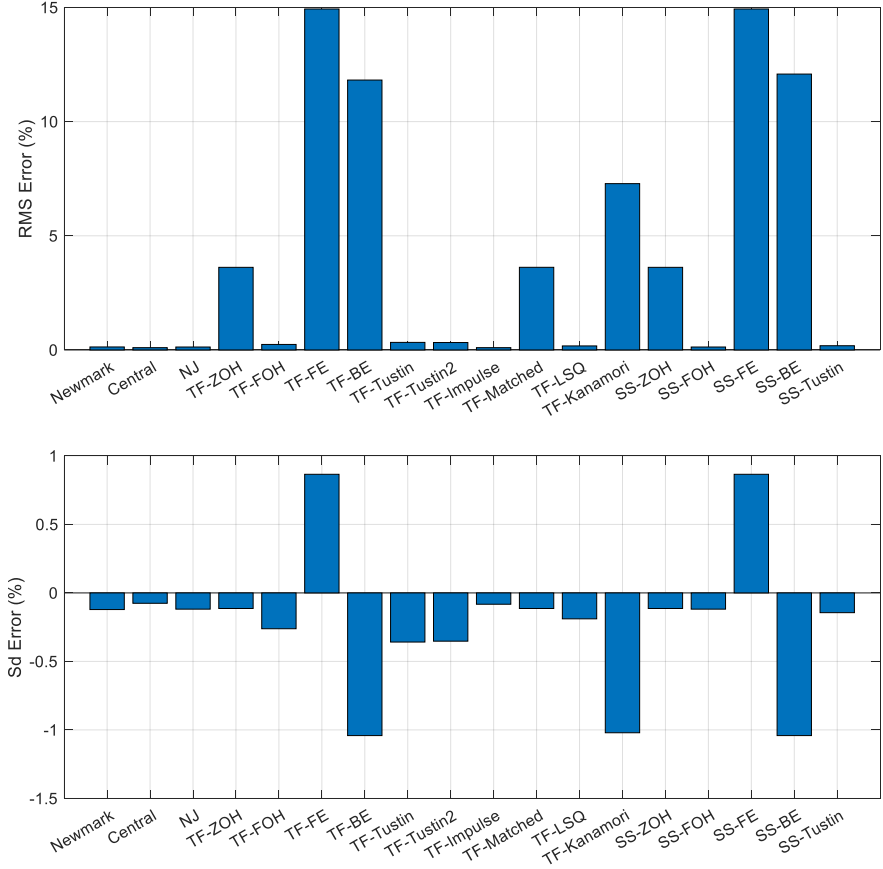

**Figure 17.** RMS and maximum displacement errors of various discretization methods ($T_n = 0.3\ sec., \xi = 5\%, \omega_0 = 5\omega_n, \Delta t_a = 0.001\ sec.$, sinc interpolation)

As a final evaluation, the time-history response (relative displacement) of a series of SDOF systems with periods listed in Table 9 and a 5% damping ratio is calculated under sinusoidal base acceleration with unit amplitude, excitation frequency equal to the natural frequency of each system, and a sampling rate of 100 Hz. The analysis is carried out for 200 seconds, and the error indices introduced in Eqs. (155) and (156) are computed for each of the 18 numerical methods. **Figure 18** presents the results of this evaluation study for three cases of analysis sampling time:

1. $\Delta t_a = \Delta t = 0.01\ sec.$,
2. $\Delta t_a = \Delta t/10 = 0.001\ sec.$ with linear interpolation of the input excitation, and
3. $\Delta t_a = \Delta t/10 = 0.001\ sec.$ with Sinc interpolation of the input excitation.

In these plots, note that the maximum limit of the vertical axes decreases from top to bottom. Moreover, methods from the same family are displayed using the same line style but different colors. As expected, Euler methods exhibit the largest errors and even stability issues for short-period structures. Other methods show good accuracy, among which the LSQ method performs particularly well because the coefficients of the transfer function are computed specifically for each structure. Excluding the LSQ method, transfer-function- or state-space-based models with FOH discretization (or its equivalent form, the Nigam–Jennings method) are among the most accurate and reliable approaches for time-history analysis.

It should also be emphasized that reducing the analysis sampling time for short-period structures is crucial for accuracy. Based on **Figure 18**, the error in both maximum displacement and RMS can be reduced to below 1% when using $\Delta t_a = 0.001\ sec.$ together with Sinc interpolation during input excitation down-sampling.

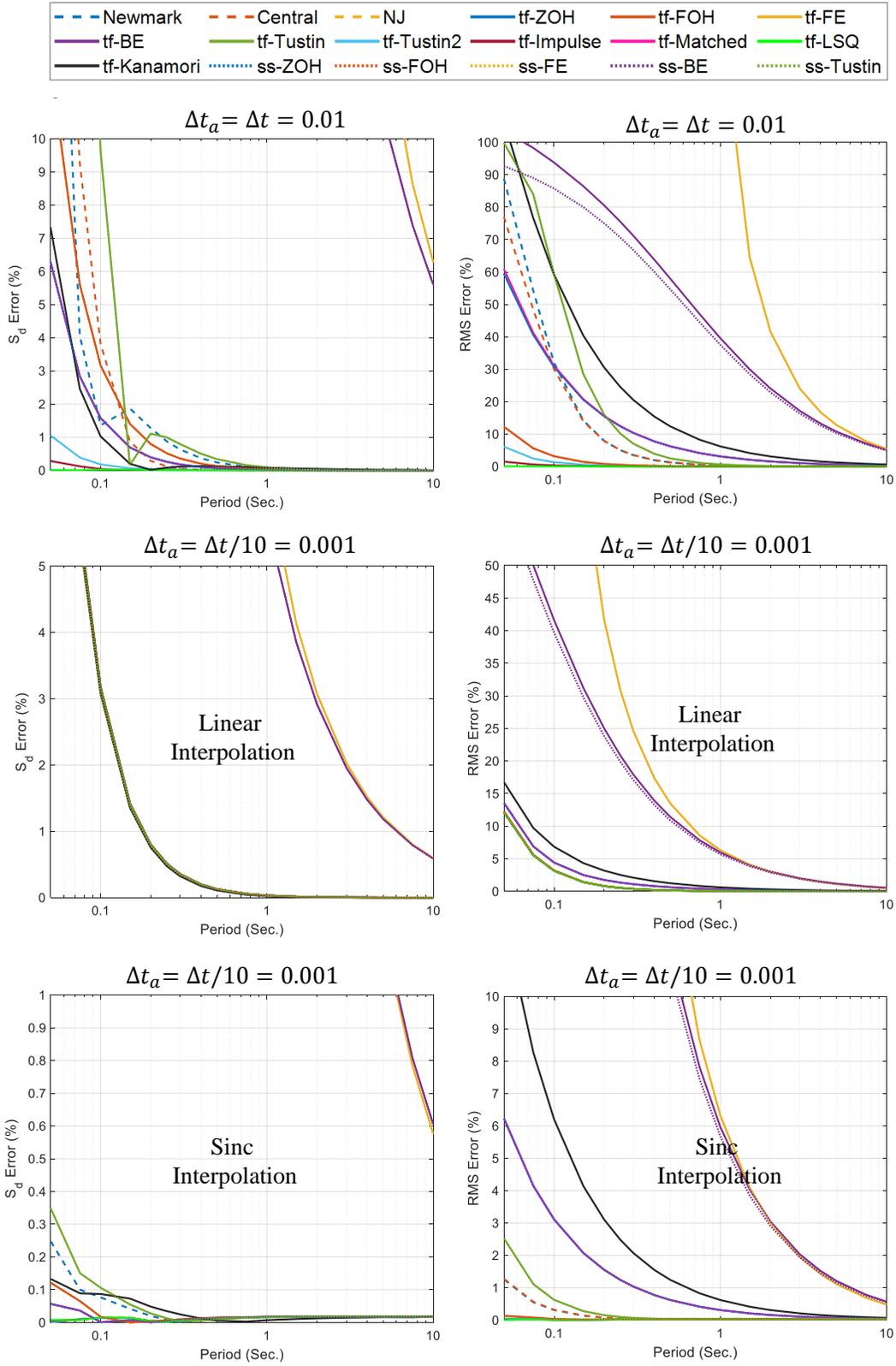

**Figure 18.** RMS and maximum displacement errors of different discretization methods for a series of structures under resonance conditions.

# CONCLUSIONS

This paper presented a unified perspective on numerical discretization methods for the seismic response analysis of linear elastic SDOF systems. By linking earthquake engineering practices to classical transfer-function- and state-space–based formulations, it was shown that well-known methods such as those of Nigam and Jennings, the Central Difference scheme, and the Newmark family were special cases of broader discretization approaches established in other engineering fields. The stability study demonstrated that most methods preserved the stability of the continuous-time system, while the comparative examples revealed that accuracy depended on both the discretization method and the analysis step size relative to the structural period. Among the examined techniques, FOH-based approaches and the Least-Squares method consistently provided accurate results, while Euler-based methods exhibited limitations. The role of input excitation processing was also highlighted, with both linear and Sinc interpolation enabling accurate comparisons with analytical solutions. Overall, the unified framework clarified the theoretical foundations of widely used numerical methods and provided guidance for their effective application in seismic response analysis and real-time implementation.

# ACKNOWLEDGMENTS

This study was partly supported by the California Geological Survey (CGS), which is gratefully acknowledged. Any opinions, findings, and conclusions expressed in this material are solely those of the author and do not reflect the views of the California Department of Conservation.